\documentclass[11pt,twoside]{article}

\usepackage{amsbsy,amsfonts,amsmath,amssymb,eucal,mathrsfs}
\usepackage[all]{xy}
\usepackage[colorlinks=true]{hyperref}

\def\no{\noindent}

\renewcommand{\ge}{\geqslant}
\renewcommand{\le}{\leqslant}

\def\itemn#1{\item[\hspace{0.6mm} {\rm (#1)}]}
\def\bin#1#2{\mbox{\scriptsize $\big\{\!\!\!\begin{array}{c} #1 \\ #2
    \end{array}\!\!\!\big\}$\normalsize}}

\def\too{\longrightarrow}

\def\into{\hookrightarrow}

\renewcommand{\tilde}{\widetilde}
\renewcommand{\bar}{\overline}

\def\lb{\lambda}
\newcommand \bbf[1]{\textbf{\textit{#1}}{}}
\newcounter{cpt}
\newcommand\tiret{\stepcounter{cpt} \bigskip \noindent (\thecpt ) }

\def\egaldef{\stackrel{{\rm df}}{=}}

\newtheorem{counter}[subsubsection]{$\!\!$}
\newtheorem{subcounter}[subsection]{$\!\!$}
\newenvironment{defi}{\begin{counter} \rm {\bf Definition.}}{\end{counter}}

\newenvironment{assu}{\begin{counter} \rm {\bf Assumption.}}{\end{counter}}

\newenvironment{prop}{\begin{counter} {\bf Proposition.}}{\end{counter}}
\newenvironment{lemm}{\begin{counter} {\bf Lemma.}}{\end{counter}}

\newenvironment{coro}{\begin{counter} {\bf Corollary.}}{\end{counter}}
\newenvironment{theo}{\begin{counter} {\bf Theorem.}}{\end{counter}}

\newenvironment{rema}{\begin{counter} \rm {\bf Remark.}}{\end{counter}}

\newenvironment{noth}{\begin{counter} \rm}{\end{counter}}
\newenvironment{proo}{{\flushleft \bf Proof:}}{\hfill $\square$ \vspace{5mm}}

\DeclareMathOperator{\Hom}{Hom}
\DeclareMathOperator{\End}{End}
\DeclareMathOperator{\Ext}{Ext}
\DeclareMathOperator{\HH}{H}

\DeclareMathOperator{\Gr}{Gr}
\DeclareMathOperator{\FG}{FG\!}
\DeclareMathOperator{\FS}{FS\!}
\DeclareMathOperator{\Id}{Id}

\DeclareMathOperator{\Spec}{Spec}

\DeclareMathOperator{\Frame}{Fr}

\DeclareMathOperator{\ffl}{fl}

\DeclareMathOperator{\Fund}{Fund}

 \def\cE{{\mathcal E}} \def\cF{{\mathcal F}}
\def\cG{{\mathcal G}}  \def\cI{{\mathcal I}}
  
  \def\cO{{\mathcal O}}

\def\cV{{\mathcal V}} \def\cW{{\mathcal W}}

\renewcommand\AA{\mathbb{A}} \newcommand\BB{\mathbb{B}}

\newcommand\GG{\mathbb{G}}

 \newcommand\NN{\mathbb{N}}
 
\newcommand\QQ{\mathbb{Q}} 
 
\newcommand\UU{\mathbb{U}} \newcommand\VV{\mathbb{V}}
\newcommand\WW{\mathbb{W}} \newcommand\XX{\mathbb{X}}
\newcommand\YY{\mathbb{Y}} \newcommand\ZZ{\mathbb{Z}}

 \def\sE{\mathscr{E}} \def\sF{\mathscr{F}}
\def\sG{\mathscr{G}}  
 \def\sK{\mathscr{K}} 
\def\sM{\mathscr{M}}  
  \def\sR{\mathscr{R}}
\def\sS{\mathscr{S}}

\begin{document}

\begin{center}
{\Large \bf Sekiguchi-Suwa theory revisited}
\end{center}

\begin{center}
Ariane M\'ezard, Matthieu Romagny, Dajano Tossici
\end{center}

\begin{center}
{\em \today}
\end{center}

\bigskip
\bigskip

\begin{center}
{\bf Abstract}

\bigskip

\begin{minipage}{16cm}
{\small We present an account of the construction by S. Sekiguchi
and N. Suwa of a cyclic isogeny of affine smooth group schemes
unifying the Kummer and Artin-Schreier-Witt isogenies. We complete
the construction over an arbitrary base ring. We extend the
statements of some results in a form adapted to a further
investigation of the models of the group schemes of roots of
unity.}
\end{minipage}
\end{center}

\bigskip

\tableofcontents

\newpage

Given a prime $p$ and an integer $n\ge 1$, consider the problem of
describing \'etale cyclic coverings of order $p^n$ of algebras, or schemes.
Over a field of characteristic $0$, the Kummer isogeny provides such a
covering which is universal on local rings. Over a field of
characteristic~$p$, an isogeny with the same virtues is given by the
Artin-Schreier-Witt theory. In the end of the nineties, T.~Sekiguchi and
N.~Suwa gave the construction of an isogeny of smooth affine
$n$-dimensional group schemes over a discrete valuation ring of mixed
characteristics, putting the Kummer isogeny and the Artin-Schreier-Witt
isogeny into a continuous family satisfying a certain universality property.
This is presented in the papers \cite{SS1} and \cite{SS2} and we give
a more detailed overview in Section~\ref{ss:over} below.

The present paper is an account of this construction, with emphasis
on some features that we found especially interesting.
We have three main goals in writing such an account.

Our first goal is to generalize their theory in such a way that it can
handle as many isogeny kernels as possible. The point is that, while Sekiguchi
and Suwa are mainly interested in {\em one} model of $\mu_{p^n}$, we are
interested in {\em all} models. For this, we need to give some
complements to the papers \cite{SS1} and \cite{SS2} and make sure that
the proofs of the generalized statements work. The result is
Theorem~\ref{th:universal_Kummer}. Also, since
the article \cite{SS2} was never published, we wanted to check thoroughly
all the details so as to rely safely on it.

Our second goal is to emphasize the geometric nature of the construction.
Indeed, the assumption that the base is a discrete valuation ring is almost
useless in \cite{SS1} and \cite{SS2}. With suitable formulations, everything
works over an (almost) arbitrary $\ZZ_{(p)}$-algebra, and the result is a
parameterization of a nice family of affine smooth group schemes called
{\em filtered group schemes}, containing plenty of models of $\mu_{p^n}$.
The parameter space is a countable union of schemes of finite type over
$\ZZ_{(p)}$, as we prove in Theorem~\ref{th:universal_alg}. We show how to
formulate things in this geometric, functorial way.

Our third goal is to propose a hopefully pleasant exposition of the theory,
with the idea that this tremendous piece of algebra deserves to be best-known.
We introduce some terminology for important concepts when we think that it
may be enlightening (fundamental morphisms, framed group schemes, Kummer
subgroup). We focus on key points rather than lengthy calculations.
We emphasize the inductive
nature of the intricate constructions with an algorithmic presentation.
We do not claim that reading our text is a gentle stroll leading without
effort to a transparent understanding of the papers \cite{SS1} and \cite{SS2}.
Rather, we hope that having a slightly different viewpoint will help the
interested reader to immerse into these papers.

\bigskip

\no {\bf Summary of contents.}
We first present the main lines of the strategy of Sekiguchi and Suwa to
construct some affine smooth group schemes embodying the unification of
Kummer and Artin-Schreier-Witt theories. (\S 1.1-1.3). Our aim is to
describe as many isogenies as possible between these groups, and to
study their kernels (\S 1.4). We recall the necessary notions on Witt vectors
(\S 2). We define and classify framed formal groups by a universal object
(Theorem 3.2.9). We emphasize that the construction by induction is
given by an explicit and computable algorithm (\S 3). Section 4 is devoted
to framed group schemes. In order to obtain algebraic objects we have to
truncate carefully the previous formal objects. At last, we consider
explicit isogenies between framed group schemes and we obtain the condition
to define finite flat Kummer group schemes~(\S 5).

\bigskip

\no {\bf Notations.} The roman and the greek alphabets do not contain
enough symbols for Sekiguchi-Suwa theory. Using the same letters for
different objects could not always be avoided. We tried our best to choose
good notations, but in some places they remain very heavy. In other places,
we changed slightly the notations of Sekiguchi and Suwa. We apologize for
the inconvenience.

\bigskip

\no {\bf Acknowledgements.} For several useful comments and conversations,
we thank Pierre Cartier, Laurent Fargues, Michel Raynaud, and Noriyuki Suwa.
We also thank
Guillaume Pagot and Michel Matignon who provided us with their notes on
the article \cite{SS2}. The first and second authors especially enjoyed
a stay in the Scuola Normale Superiore di Pisa where part of this work
was done. The third author had fruitful stays at the MPIM in Bonn, at the
IHES in Bures-sur-Yvette, and spent some time in Paris to work on this
project invited by the University Paris 6, the University of
Versailles Saint-Quentin and the IHP, during the Galois Trimester.
The three authors also spent a very nice week in the CIRM in Luminy.
We thank all these institutions for their support and hospitality.

\section{Overview of Sekiguchi-Suwa theory} \label{ss:over}

\subsection{Unifying Kummer and Artin-Schreier-Witt theories}

Fix a discrete valuation ring $R$ with fraction field $K$ of
characteristic $0$ and residue field $k$ of characteristic $p>0$.
Let $W_n$ be the scheme of Witt vectors of length $n$ and $\GG_m$ the
multiplicative group scheme. The work of Sekiguchi and Suwa provides an explicit
construction of an isogeny $\cW_n\too\cV_n$ of smooth affine $n$-dimensional
group schemes over $R$ with special fibre isomorphic to the
Artin-Schreier-Witt isogeny
$$
\wp:W_{n,k}\too W_{n,k} \ , \ x\mapsto x^p-x,
$$
and generic fibre isomorphic to the Kummer-type isogeny
$$
\Theta:(\GG_{m,K})^n\too (\GG_{m,K})^n \ , \
(x_1,\dots,x_n) \mapsto (x_1^p,x_2^px_1^{-1},\dots,x_n^px_{n-1}^{-1}),
$$
such that any $p^n$-cyclic finite \'etale extension of local flat $R$-algebras
is obtained by base change from $\cW_n\too\cV_n$. The isogeny $\Theta$
is essentially equivalent to the usual one-dimensional isogeny
$x\mapsto x^{p^n}$ for the purposes of Kummer theory, and is of course
best-suited to the unification with the Artin-Schreier-Witt theory.
In the strategy of Sekiguchi and Suwa to complete this goal,
let us single out three steps:
\begin{trivlist}
\itemn{A} Describe a family of smooth $n$-dimensional group schemes that are good
candidates to be the domain and target of the sought-for isogeny
(this is done in Sections 3, 4, 5 of \cite{SS2}). These are called
{\em filtered group schemes}.
\itemn{B} Choose suitably the parameters in the previous constructions
so as to produce a group scheme $\cW_ n$ (Section~8 of \cite{SS2}) with a
finite flat subgroup scheme $(\ZZ/p^n\ZZ)_R$. This step requires $R$
to contain the $p^n$-th roots of unity.
\itemn{C} Compute the group $\cV_n=\cW_n/(\ZZ/p^n\ZZ)$ and the isogeny
$\cW_n\to\cV_n$ (Section~9 of \cite{SS2}).
\end{trivlist}
We will now present these steps in a little more detail.


\subsection{Filtered group schemes}

Let us briefly describe Step~(A), the description of the family of smooth
group schemes relevant to the problem.
The groups are constructed with two guiding principles: firstly they are models
of $(\GG_{m,K})^n$, and secondly they are extensions of a group of the same
type in dimension one less by a $1$-dimensional group, in the same way as
$W_{n,k}$ is an extension of $W_{n-1,k}$ by $\GG_{a,k}$. For $n=1$, the
smooth models of $\GG_{m,K}$ with connected fibres are known as some group
schemes $\cG^{\lambda}=\Spec(R[X,1/(1+\lambda X)])$, where $\lambda\in R$
is a parameter (see the papers~\cite{WW} and \cite{SOS}). Thus we are
led to consider {\em filtered group schemes of type $(\lambda_1,\dots,\lambda_n)$}
for various $n$-tuples of elements $\lambda_i\in R$, defined recursively
as the extensions of a group $\cE$ of type $(\lambda_1,\dots,\lambda_{n-1})$
by the group~$\cG^{\lambda_n}$. We see that in order to obtain the
$n$-dimensional group schemes, we have to describe the group
$\Ext^1(\cE,\cG^{\lambda})$ classifying such extensions. This is easy
when $\lambda$ is invertible i.e. $\cG^{\lambda}\simeq\GG_{m,R}$,
since one can prove easily by d\'evissage that $\Ext^1(\cE,\GG_{m,R})=0$.
Therefore, in order to understand $\Ext^1(\cE,\cG^{\lambda})$ we must
measure the difference between $\cG^{\lambda}$ and $\GG_{m,R}$. This
is done with an exact sequence of sheaves on the small flat site
$$
0 \too \cG^{\lambda}\too \GG_{m,R}
\stackrel{\rho}{\too} i_*\GG_{m,R/\lambda}\too 0
$$
where $i:\Spec(R/\lambda R)\too\Spec(R)$ is the closed immersion
(we make the convention that all sheaves supported on the empty set
are $0$, e.g. $i_*\GG_{m,R/\lambda}=0$ if $\lambda$ is invertible).
The long exact sequence for the functor $\Hom(\cE,\cdot)$ gives
$$
\Ext^1(\cE,\cG^{\lambda})\simeq
\Hom(\cE_{R/\lambda},\GG_{m,R/\lambda})/\rho_*\Hom(\cE,\GG_{m,R})
$$
(and $\Hom(\cE_{R/\lambda},\GG_{m,R/\lambda})=0$ if $\lambda$ is
invertible, according to our previous convention).

At this point, the problem becomes essentially to describe
$\Hom(\cE_{R/\lambda},\GG_{m,R/\lambda})$. Technically, this is one of
the key points
of Sekiguchi and Suwa's work. This group of homomorphisms is parameterized
by a suitable generalization of the classical Artin-Hasse exponential
series. It is therefore really in the formal world that the crucial objects
live, as formal power series satisfying the important identities. Accordingly,
the formal theory (the construction of filtered formal groups) precedes,
and is the inspiration for, the algebraic theory (the construction of
filtered group schemes). Here, it is worth pointing out that the
construction of extensions in the formal case takes a slightly different
turn, because no analogue of the exact sequence
$0\to\cG^\lambda\to\GG_m\to i_*\GG_m\to 0$
is available. Instead one considers the composition
$$
\partial:\Hom(\hat\cE,\hat\GG_m)\too H^2_0(\hat\cE,\hat\cG^\lambda)
\too \Ext^1(\hat\cE,\hat\cG^\lambda)
$$
that associates to a morphism a Hochschild $2$-cocyle and then the
extension it gives birth to. The point is that in the algebraic case, the map
$\partial$ is obtained as the coboundary of a long exact cohomology
sequence which is {\em not} available in the formal case, while in the formal
case the map $\partial$ is obtained using the {\em surjective} map
$H^2_0(\hat\cE,\hat\cG^\lambda)\to \Ext^1(\hat\cE,\hat\cG^\lambda)$ which
tends to be {\em zero} in the algebraic case.

\subsection{Finite flat subgroup schemes} \label{ss:ff_gp_sch}

Let us now make some comments on Steps (B) and (C). Filtered group schemes
$\cE$ have {\em filtered subgroup schemes}, obtained by successive extensions
of subgroups. We will see that their construction provides natural morphisms
$\alpha:\cE\to (\GG_m)^n$ that are {\em model maps}, that is to say,
isomorphisms on the generic fibre. On the generic fibre, these morphisms provide
natural filtered subgroup schemes of $\cE_K$ isomorphic to $\mu_{p^n,K}$: one
just has to pullback via $\alpha$ the kernel of the Kummer isogeny
$\Theta_K:(\GG_{m,K})^n\to (\GG_{m,K})^n$.
By taking the closure in $\cE$, one produces interesting candidates to be
finite flat models of $\mu_{p^n,K}$. If $R$ contains the
$p^n$-roots of unity, and for suitable choices of the parameters of the
extensions, one obtains a filtered group scheme $\cE=\cW_n$ and a model
of $\mu_{p^n,K}\simeq (\ZZ/p^n\ZZ)_K$ which turns out to be the constant
group $(\ZZ/p^n\ZZ)_R$. Sekiguchi and Suwa specialize to this case and
study the quotient isogeny. They prove that these
objects realize the unification of the Kummer and Artin-Schreier-Witt
exact sequences.

\subsection{Our presentation of the theory}

Our personal interest does not lie in {\em one} single model of $\mu_{p^n,K}$
but in {\em all} possible models one can exhibit (see the article~\cite{MRT}).
It is therefore very important for us to leave the parameters as free as possible.
We call {\em Kummer subschemes} the subschemes $G$ obtained by scheme-theoretic
closure in the way described in~\ref{ss:ff_gp_sch}. Then the framework of
Sekiguchi and Suwa allows to characterize when a Kummer subscheme is
finite locally free over the base ring $R$. In fact, the 'good' object is the
isogeny $\cE\to \cF=\cE/G$ itself, and we are able to construct isogenies
between filtered group schemes, whose kernels are the finite flat models of
$\mu_{p^n,K}$ we are interested in.

If we incorporate the various choices of parameters into the definitions, we
obtain a notion of {\em framed group scheme} whose moduli problem is
(tautologically) representable by a scheme. This scheme is a nice parameter
space for filtered group schemes. It has a formal and an algebraic version.
We formulate things with this vocabulary.

Finally, we point out that almost no restriction on the base ring $R$ is
necessary. In particular, it need not be a discrete valuation ring, not even
an integral domain. The only important point is that the parameters
$\lambda_i$ of the successive extensions should be nonzerodivisors. Thus we
work throughout with an arbitrary $\ZZ_{(p)}$-algebra.

\section{Witt vectors} \label{ss:WittV}

The prime number $p$ is fixed. This section is devoted to generalities on the
ring scheme of Witt vectors $W$. We first recall basic notations concerning $W$
and some of its endomorphisms. Then we define the formal completion of $W$
and study its stability under the endomorphisms defined before. Finally we
introduce various objects related to the scheme of Witt vectors over the affine line.
As a general rule, we keep the notations of the papers
\cite{SS2} and \cite{SS1}.

\subsection{Witt vectors}

We briefly indicate our notations for the ring scheme $W$ of Witt vectors over
the integers. The letters $X,Y,Z,A$ denote infinite vectors of
indeterminates, with $X=(X_0,X_1,\dots)$, etc.

\begin{noth} {\bf Ring scheme structure.}
The {\em scheme of Witt vectors} is $W=\Spec(\ZZ[Z_0,Z_1,\dots])$. Its structure
is defined using the {\em Witt polynomials} defined for all integers $r\ge 0$
by:
$$
\Phi_r(Z)=\Phi_r(Z_0,\dots,Z_r)=Z_0^{p^r}+pZ_1^{p^{r-1}}+\dots+p^rZ_r.
$$
The addition and multiplication of the Witt ring scheme are defined
respectively, on the function ring level, by the assignments
$Z_r\mapsto S_r(X,Y)$ and $Z_r\mapsto P_r(X,Y)$,
where
$$
S_r(X,Y)=S_r(X_0,\dots,X_r,Y_0,\dots,Y_r)
\ ,\
P_r(X,Y)=P_r(X_0,\dots,X_r,Y_0,\dots,Y_r)
$$
are the unique polynomials with integer coefficients satisfying
for all $r\ge 0$ the identities:
$$
\begin{array}{l}
\Phi_r(S_0(X,Y),\dots,S_r(X,Y))=\Phi_r(X_0,\dots,X_r)+\Phi_r(X_0,\dots,Y_r), \medskip \\
\Phi_r(P_0(X,Y),\dots,P_r(X,Y))=\Phi_r(X_0,\dots,X_r)\,\Phi_r(X_0,\dots,Y_r). \\
\end{array}
$$
\end{noth}

\begin{noth} {\bf Frobenius, Verschiebung, Teichm\"uller, $T$ map.}
The ring scheme endomorphism $F:W\to W$ called {\em Frobenius} is defined
by the assignment $Z_r\mapsto F_r(X)$, where the
$F_r(X)=F_r(X_0,\dots,X_{r+1})$ are the unique polynomials satisfying
for all $r\ge 0$ the identities:
$$
\Phi_r(F_0(X),\dots,F_r(X))=\Phi_{r+1}(X_0,\dots,X_{r+1}).
$$
The additive group scheme endomorphism
$$
V:W=\Spec(\ZZ[X_0,X_1,\dots])\to W=\Spec(\ZZ[Z_0,Z_1,\dots])
$$
called {\em Verschiebung}
is defined by the assignments $Z_0\mapsto 0$ and $Z_r\mapsto X_{r-1}$ for
$r\ge 1$. Let $\AA^1=\Spec(\ZZ[X_0])$ be the affine line over $\ZZ$.
Then the multiplicative morphism $[\,\cdot\,]:\AA^1\to W$ called
{\em Teichm\"uller representative} is
defined by the assignments $Z_0\mapsto X_0$ and $Z_r\mapsto 0$ if $r\ge 1$.

An important role in Sekiguchi-Suwa theory is played by the morphism
$T:W\times W\to W$ called (by us) the {\em $T$ map}, defined by the assignment
$Z_r\mapsto T_r(Y,X)$, where the $T_r(Y,X)$ are the unique polynomials
satisfying for all $r\ge 0$ the identities:
$$
\Phi_r(T_0(Y,X),\dots,T_r(Y,X))=Y_0^{p^r}\Phi_r(X)
+pY_1^{p^{r-1}}\Phi_{r-1}(X)+\dots+p^rY_r\Phi_0(X).
$$
Existence and uniqueness of the sequence $T(Y,X)=(T_0(Y,X),T_1(Y,X),\dots)$
are granted by Bourbaki \cite{B}, \S~1, no.~2, Prop.~2, applied to the ring
$\ZZ[Y,X]$ endowed with the endomorphism $\sigma$ raising each variable to
the $p$-th power. Note that in \cite{SS2} the
notation for $T(Y,X)$ is $T_YX$, a notation that we will also use.
The morphism $T$ is additive in the
second variable i.e. gives rise to a morphism $T:W\to\End(W,+)$.

Some of these definitions are really more pleasant in terms of functors
of points. This is typically the case for the morphisms $V$, $T$ and
$[\,\cdot\,]$. Let us indicate them: given a ring $A$ and Witt vectors
$a,x\in W(A)$, we have $V(x)=(0,x_0,x_1,\dots)$, $[x_0]=(x_0,0,0,\dots)$ and
$T_ax=\sum_{r\ge 0}\, V^r([a_r]x)$, see \cite{SS2}, Lemma 4.2.
\end{noth}

\subsection{Formal completion}

The formal completion of the group scheme of Witt vectors along the zero
section is the subfunctor $\hat W\subset W$ defined by:
$$
\hat W(A)\egaldef\big\{\bbf{a}=(a_0,a_1,a_2,\dots)\in W(A),
\mbox{ $a_i$ nilpotent for all $i$},
\mbox{ $a_i=0$ for $i\gg 0$}\big\}.
$$
This is the completion in Cartier's sense (see~\cite{Ca});
note that in infinite dimension, several reasonable definitions of
completion exist (a different one may be found for example in~\cite{Ya},
example 3.24).

\begin{lemm} \label{lm:ideal}
The formal completion $\hat{W}$ is an ideal of $W$.
\end{lemm}

\begin{proo}
We introduce a filtration of $\hat W$ by subfunctors $\hat W_{M,N}$
($M,N\ge 1$ integers) with
$$
\hat W_{M,N}(A)=\big\{\bbf{a}\in W(A),\,a_i=0 \mbox{ for } i\ge M
\mbox{ and } (a_i)^N=0 \mbox{ for } i\ge 0\big\}.
$$
It is clear that this filtration is exhaustive. Hence it is enough to prove
that for all $M,N$ there exist $M',N'$ such that
$\hat W_{M,N}+\hat W_{M,N}\subset \hat W_{M',N'}$
and $W\times \hat W_{M,N}\subset \hat W_{M',N'}$. The proof in the two cases
is very similar, so we will treat only the case of the sum.

\smallskip

\no {\em Step 1: we may assume that $p$ is invertible in the base ring.}
Indeed, $\hat W_{M,N}$ is a closed subfunctor of $\hat{W}$ which is representable
by a finite flat $\ZZ$-scheme. So if the addition map on $\hat W_{M,N}$ factors
over $\ZZ[1/p]$ through some $\hat W_{M',N'}$, then by taking scheme-theoretic
closures one finds that it factors through $\hat W_{M',N'}$ over $\ZZ$ as well.

\smallskip

\no {\em Step 2: let $(X,Y)$ be the universal point of $W_{M,N}\times W_{M,N}$,
where $X=(X_0,X_1,\dots)$ and $Y=(Y_0,Y_1,\dots)$. Then the coefficients of
the sum $S=X+Y$ are nilpotent.} This is clear, since for each $i$ the
coefficient $S_i$ is a polynomial in the $X_j,Y_j$.

\smallskip

\no {\em Step 3: for all $i\ge r:=M-1+\log_p(N)$, we have
$\Phi_i(S)=0$.} Indeed, we have
$$
\Phi_i(X)=\sum_{j=0}^i\,p^j(X_j)^{p^{i-j}}=\sum_{j=0}^{M-1}\,p^j(X_j)^{p^{i-j}}=0
$$
since $j\le M-1$ implies that $p^{i-j}\ge p^{i-M+1}\ge p^{r-M+1}\ge N$.
Similarly we have $\Phi_i(Y)=0$ and hence $\Phi_i(S)=\Phi_i(X)+\Phi_i(Y)=0$.

\smallskip

\no {\em Step 4: by Step 2, let $P$ be such that $(S_0)^P=\dots=(S_{r-1})^P=0$.
Then $S_i=0$ for all $i\ge \log_p((P-1)(1+p+\dots+p^{r-1}))$.}
For the weight of Witt vectors giving weight $p^i$ to $X_i$ and $Y_i$, the
element $S_i$ is
homogeneous of weight $p^i$. Since $p$ is invertible, using Step 3 and
induction we see that for all $i\ge r$, the element $S_i$ is a polynomial
in $S_0,\dots,S_{r-1}$. By the choice of $P$, a monomial
$(S_0)^{j_0}\dots (S_{r-1})^{j_{r-1}}$ will be nonzero only if all exponents
$j_0,\dots,j_{r-1}$ are less than $P-1$, hence the weight is
$j_0+pj_1+\dots+p^{r-1}j_{r-1}\le (P-1)(1+p+\dots+p^{r-1})$.
We get the claim by contraposition.

\smallskip

\no {\em Step 5: conclusion.} By Step 4, we can take
$M'=\log_p((P-1)(1+p+\dots+p^{r-1}))$ and the existence of $N'$ is given by Step 2.
\end{proo}

\begin{rema}
In Sections~\ref{sec:AT} and \ref{sec:KS}, we try to give a presentation of
Sekiguchi-Suwa theory adapted to computations. In particular, in
Lemma~\ref{lm:TDE} we give an explicit degree of truncation for the Artin-Hasse
exponentials that is sufficient to compute filtered group schemes. It is equally
desirable to have explicit bounds for the number of nonzero terms of the Witt
vectors that appear, but this desire is in fact limited by the difficulty to
give a reasonably explicit bound for the number of nonzero coefficients of
the {\em sum} of two Witt vectors, as we saw in the proof of Lemma~\ref{lm:ideal}.
\end{rema}

\begin{lemm} \label{lm:F_stable}
The formal completion $\hat{W}$ is stable under $F$ and $V$.
\end{lemm}

\begin{proo}
For $V$ there is nothing to say, and for $F$ the strategy of the proof
of Lemma~\ref{lm:ideal} works almost unchanged.
\end{proo}

\begin{lemm} \label{lm:T_stable}
Let $W^f$ be the subfunctor of $W$ composed of Witt vectors with finitely
many nonzero coefficients. Then $T$ induces a morphism
$W^f\times\hat{W}\to \hat{W}$.
\end{lemm}

We point out that $W^f$ has no (additive or whatever) structure.

\begin{proo}
Using the formulas $T_ax=\sum_{r\ge 0}\, V^r([a_r]x)$ and
$[a]x=(ax_0,a^px_1,a^{p^2}x_2,\dots)$, this is obvious.
\end{proo}

\subsection{Witt vectors over the affine line} \label{ss:W_A1}

Let $\AA^1=\Spec(\ZZ[\Lambda])$ be the affine line over the integers, and let
$i:\Spec(\ZZ)\into \AA^1$ be the closed immersion of the origin, given by
$\Lambda=0$. In the paper \cite{SOS}, the study of the multiplicative group
scheme over the affine line leads to introduce a certain group scheme
$\cG^\Lambda$ (the notation in {\em loc. cit.} is $\cG^{(\lambda)}$).
In this section, we expand the idea behind the introduction of this group scheme,
because we notice that when we consider a group scheme over $\AA^1$
(favourite examples are $\GG_m$ or $W$), the groups of elements
{\em vanishing at the origin} and those {\em supported at the origin}
are especially important. In this way, we introduce a $W$-module scheme
$W^\Lambda$. We recall the definition of $\cG^\Lambda$ which fits in the same
framework. Note that we simplify the notations $F^{(\Lambda)}$,
$\cG^{(\Lambda)}$, $\alpha^{(\Lambda)}$ from the papers \cite{SOS} and
\cite{SS2} to $F^\Lambda$, $\cG^\Lambda$, $\alpha^\Lambda$.

\begin{prop} \label{pp:cG_cW}
Let $\AA^1=\Spec(\ZZ[\Lambda])$ be the affine line over the integers, and let
$i:\Spec(\ZZ)\into \AA^1$ be the closed immersion given by $\Lambda=0$.
Let $\AA^1_{\ffl}$ denote the small flat site of $\AA^1$.
\begin{trivlist}
\itemn{1} The canonical morphism $\GG_m\to i_*\GG_m$ fits into an exact
sequence
$$
0\too\cG^\Lambda\stackrel{\alpha^\Lambda}{\too}\GG_m\too i_*\GG_m\too 0
$$
of abelian sheaves on $\AA^1_{\ffl}$,
where $\cG^\Lambda$ is a flat commutative group scheme.
\itemn{2} The canonical morphism $W\to i_*W$ fits into
an exact sequence
$$
0\too W^\Lambda\stackrel{\Lambda}{\too} W\too i_*W\too 0
$$
of abelian sheaves on $\AA^1_{\ffl}$, where $W^\Lambda$ is a flat
$W$-module scheme. Here, the scheme $W^\Lambda$ has the same underlying
scheme as $W$ and the first map is
$$
x=(x_0,x_1,x_2,\dots)\mapsto
\Lambda.x:=(\Lambda x_0,\Lambda x_1,\Lambda x_2,\dots).
$$
\end{trivlist}
\end{prop}

An algebra $R$ and an element $\lambda\in R$ define an $R$-point
$\Spec(R)\to\AA^1$. The pullbacks of $\alpha^\Lambda$ and
$\Lambda:W^\Lambda\to W$ along this point give a morphism of $R$-group
schemes which we will denote $\alpha^{\lambda}:\cG^{\lambda}\to\GG_m$
and a morphism of $R$-schemes in $W$-modules which we will denote
$\lambda:W^{\lambda}\to W$.

\begin{proo}
We treat only case (2), since case (1) is similar and even simpler.
The scheme $W^\Lambda$ and the map $\Lambda:W^\Lambda\to W$ are defined
in the statement. These fit into an exact sequence, functorial in the flat
$\ZZ[\Lambda]$-algebra $R$:
$$
0\too W^\Lambda(R)\stackrel{\Lambda}{\too} W(R)
\too (i_*W)(R)=W(R/\Lambda R)\too 0.
$$
Thus the map $\Lambda$ identifies $W^\Lambda(R)$ with the ideal of
$W(R)$ of vectors all whose components are multiples of~$\Lambda$.
It follows that for all $u,v\in W^\Lambda(R)$ and $a\in W(R)$, the sum
$u+v$ and the product $au$, computed in $W(R)$, again lie in this ideal. By
taking for $R$ the function ring of~$W^\Lambda$, we see that the universal
polynomials giving Witt vector addition and multiplication
$$
S_0(\Lambda.u,\Lambda.v),S_1(\Lambda.u,\Lambda.v),S_2(\Lambda.u,\Lambda.v),\dots
\quad,\quad
P_0(\Lambda.a,\Lambda.v),P_1(\Lambda.a,\Lambda.v),P_2(\Lambda.a,\Lambda.v),\dots
$$
are divisible by $\Lambda$, that is $S_i(\Lambda.u,\Lambda.v)=\Lambda S'_i(u,v)$
and $P_i(\Lambda.a,\Lambda.u)=\Lambda P'_i(a,u)$. By flatness, the polynomials
$S'_i$ and $P'_i$ are uniquely determined and they define the $W$-module
structure on the scheme $W^\Lambda$.
\end{proo}

\begin{rema}
We could also define $W^\Lambda$ and $\cG^\Lambda$ as dilatations of $W$ and
$\GG_m$ along the respective unit sections of the special fibre $\Lambda=0$.
When the base ring is a discrete valuation ring $R$, the dilatation of an
$R$-scheme $X$ along a closed subscheme of the special fibre is defined in
Chapter~3 of \cite{BLR}. The same construction works in the following more
general setting. Consider a base scheme $S$, a Cartier divisor $S_0=V(\cI)$,
an $S$-scheme $X$, and a closed subscheme $Y_0$ of $X_0=X\times_S S_0$. Then
there exists a morphism of $S$-schemes $u:X'\to X$ where $X'$ is an
$S$-scheme without $\cI$-torsion such that $u(X'_0)\subset Y_0$, and which
is universal with these properties. The scheme $X'$ is called the
{\em dilatation of $X$ along $Y_0$}.
\end{rema}

We close the section with a lemma that plays a key role in the
development of the theory.

\begin{lemm} \label{lm:F_flat}
Let $W$ be the ring scheme of Witt vectors over the affine line
$\AA^1=\Spec(\ZZ[\Lambda])$. Then, the additive endomorphism
$F^\Lambda:=F-[\Lambda^{p-1}]:W\to W$ is faithfully flat.
\end{lemm}

Of course, here again, for an algebra $R$ and an element $\lambda\in R$
we obtain a faithfully flat endomorphism $F^\lambda:W_R\to W_R$.

\begin{proo}
See \cite{SS1}, Proposition~1.6 and Corollaries 1.7-1.8, and \cite{SS2}, Lemma~4.5.
\end{proo}

\section{Formal theory} \label{sec:FT}

In Subsection~\ref{ss:DAHE}, we introduce the deformed Artin-Hasse
exponentials studied by Sekiguchi and Suwa. These power series satisfy
important identities that allow to construct formal filtered group schemes
by successive extensions. This is explained in~\ref{ss:fil_fgs}, with
Theorem~\ref{th:universal_formal} summarizing the main properties
of the construction.

\subsection{Deformed Artin-Hasse exponentials} \label{ss:DAHE}

In order to describe the homomorphisms from formal filtered group schemes
(introduced in Subsection~\ref{ss:fil_gs}) to the formal multiplicative
group $\hat\GG_m$, we will need some deformations of Artin-Hasse exponentials.
For simplicity, we will call them {\em deformed exponentials}. In the
non-formal case, we will also need some truncations of these series.
We introduce all these objects here.

Given indeterminates $\Lambda$, $U$ and $T$, we define a formal power series
in $T$ with coefficients in $\QQ[\Lambda,U]$ by
$$
E_p(U,\Lambda,T)=
(1+\Lambda T)^{\frac{U}{\Lambda}}\prod_{k=1}^\infty\,
(1+\Lambda^{p^k}T^{p^k})^{\frac{1}{p^k}\left(\left(\frac{U}{\Lambda}\right)^{p^k}
-\left(\frac{U}{\Lambda}\right)^{p^{k-1}}\right)} \ .
$$
It satisfies basic properties such as $E_p(0,\Lambda,T)=1$ and
$E_p(MU,M\Lambda,T)=E_p(U,\Lambda,MT)$, where $M$ is another indeterminate.
It is a deformation of the classical Artin-Hasse exponential
$E_p(T)=\prod_{k=0}^\infty\,\exp(T^{p^k}/p^k)$ in the sense that $E_p(1,0,T)=E_p(T)$.
Given a vector of indeterminates $\UU=(U_0,U_1,\dots)$, we define
a power series in $T$ with coefficients in $\QQ[\Lambda,U_0,U_1,\dots]$ by
\begin{equation}\label{eq:E_p(U,lambda,T)}
E_p(\UU,\Lambda,T)=\prod_{\ell=0}^\infty\,E_p(U_\ell,\Lambda^{p^\ell},T^{p^\ell}).
\end{equation}
It is proven in~\cite{SS1}, Corollary~2.5 that the series
$E_p(U,\Lambda,T)$ and $E_p(\UU,\Lambda,T)$ are integral at $p$,
that is, they have their coefficients in $\ZZ_{(p)}[\Lambda,U]$
and $\ZZ_{(p)}[\Lambda,U_0,U_1,\dots]$ respectively. It follows
that given a $\ZZ_{(p)}$-algebra $A$, elements $\lambda,a\in A$
and $\bbf{a}=(a_0,a_1,\dots)\in A^\NN$, we have specializations
$E_p(a,\lambda,T)$ and $E_p(\bbf{a},\lambda,T)$ which are power
series in $T$ with coefficients in $A$. We usually consider
$\bbf{a}$ as a Witt vector, i.e. as an element in $W(A)$. One must
however be aware that since $W(A)$ has the extra structure of a
ring, this introduces the slight ambiguity that
$E_p(\bbf{a},\lambda,T)$ might be interpreted as the result of
specializing $U$ to $\bbf{a}$ in the series $E_p(U,\Lambda,T)$,
resulting in a series with coefficients in $W(A)$ (note that if
$A$ is a $\ZZ_{(p)}$-algebra then so is $W(A)$). However, in
Sekiguchi-Suwa theory the symbol $E_p(\bbf{a},\lambda,T)$ always
denotes a specialization of $E_p(\UU,\Lambda,T)$ so that no
confusion can come up.

\smallskip

Now we borrow some terminology from Fourier analysis.

\begin{defi} \label{df:fund_1}
Let $A$ be a $\ZZ_{(p)}$-algebra, $\lambda\in A$ an element and $k\ge 1$ a
prime-to-$p$ integer. A series of the form $E_p(\bbf{a},\lambda,T^k)$ is
called a {\em $k$-th harmonic} and a $1$-st harmonic is also called
a {\em fundamental}. A morphism $\hat\cG^\lambda\to\hat\GG_m$ defined by a
fundamental is called a {\em fundamental morphism}.
\end{defi}

The significance of this terminology is explained by the following easy lemma.

\begin{lemm}
Let $A$ be a $\ZZ_{(p)}$-algebra and $\lambda\in A$. Then every formal
power series $G\in A[[T]]$ such that $G(0)=1$ may be decomposed uniquely
as a product of harmonics. More precisely, there exist unique vectors
$\bbf{a}_k=(a_{k0},a_{k1},\dots)\in W(A)$ for all prime-to-$p$ integers $k$,
such that $G(T)=\prod_{p\nmid k}\,E_p(\bbf{a}_k,\lambda,T^k)$.
\end{lemm}

\begin{proo}
(See Remark~2.10 of \cite{SS1}.) The claim will follow simply from the
fact that $E_p(U,\Lambda,T)\equiv 1+UT \mod T^2$.
Write $G(T)=1+g_1T+g_2T^2+\dots$ and let $v:\NN\setminus\{0\}\to\NN$ be the
$p$-adic valuation. We prove by induction on $n\ge 1$ that there exist unique
elements $b_1,\dots,b_n$ in $A$ such that
$$
G(T)E_p(b_1,\lambda^{p^{v(1)}},T)^{-1}E_p(b_2,\lambda^{p^{v(2)}},T^2)^{-1}\dots
E_p(b_n,\lambda^{p^{v(n)}},T^n)^{-1}\equiv 1 \mod T^{n+1}.
$$
For $n=1$ we have $G(T)\equiv 1+g_1T \mod T^2$ and then it is necessary and
sufficient to put $b_1=g_1$. If the claim is proven for $n\ge 1$, then we
have
$$
G(T)\prod_{i=1}^n\,E_p(b_i,\lambda^{v(i)},T^i)^{-1}\equiv 1+c_{n+1}T^{n+1}
\mod T^{n+2}
$$
for some $c_{n+1}\in A$, and it is necessary and sufficient to put
$b_{n+1}=c_{n+1}$. Finally we obtain
$$
G(T)=\prod_{i=1}^\infty\,E_p(b_i,\lambda^{p^{v(i)}},T^i)
$$
and the claim follows by defining $\bbf{a}_k:=(b_k,b_{kp},b_{kp^2},\dots)$.
\end{proo}

Let $\AA^1=\Spec(\ZZ_{(p)}[\Lambda])$ be the affine line over the
$p$-integers. We finally remark that, generalizing what happens for the
classical rtin-Hasse exponential (see \cite{SS1}, Corollary 2.9.1),
the exponential $E_p(\UU,\Lambda,T)$ gives a homomorphism
$$
W_{\AA^1}\too {\mathbf{\Lambda}}_{\AA^1},
$$
where
${\mathbf{\Lambda}}_{\AA^1}=\Spec(\ZZ_{(p)}[\Lambda,X_1,\dots,X_n,\dots])$
is the $\AA^1$-group scheme  whose group of $R$-points, for
any $\ZZ_{(p)}[\Lambda]$-algebra $R$, is the abelian
multiplicative group $1+TR[[T]]$. The above homomorphism is in fact a
closed immersion, and by the above lemma there is  an isomorphism
$$
\prod_{p\nmid k} W_{\AA^1}\stackrel{\sim}{\too}\mathbf{\Lambda}_{\AA^1}.
$$


\subsection{Construction of framed formal groups} \label{ss:fil_fgs}

Let $R$ be a $\ZZ_{(p)}$-algebra and let $\lambda_1,\lambda_2,\dots$ be
elements of $R$.

\begin{defi} \label{df:FFGS}
A {\em filtered formal $R$-group of type $(\lambda_1,\dots,\lambda_n)$}
is a sequence
$$\hat\cE_0=0,\hat\cE_1,\dots,\hat\cE_n$$
of affine smooth commutative formal group schemes such that for each
$i=1,\dots,n$ the formal group $\hat\cE_i$ is an extension of
$\hat\cE_{i-1}$ by $\hat\cG^{\lambda_i}$.
\end{defi}

We now indicate a procedure due to Sekiguchi and Suwa for constructing
filtered formal groups. It works under the following:

\begin{assu}
The elements $\lambda_1,\lambda_2,\dots$ are not zero divisors in $R$.
\end{assu}

The procedure involves some choices which we take into account by
introducing notions of {\em frames} and {\em framed formal groups}.
In this way, the refined procedure becomes universal. We adapt the
construction of~\cite{SS2} accordingly.

\smallskip

Let $W$ be the $R$-group scheme of infinite Witt vectors. For each
$\lambda\in R$, we have the morphisms of $R$-group schemes
$\alpha^{\lambda}:\cG^{\lambda}\to\GG_m$ and
${\lambda}:W^{\lambda}\to W$ introduced in
Subsection~\ref{ss:W_A1}. For each integer $n\ge 1$, we have a
product morphism
${\lambda}\times\dots\times{\lambda}:(W^{\lambda})^n\to
W^n$ which by abuse we again denote by the symbol
${\lambda}$.

\smallskip

\begin{noth} {\bf Description of the procedure.}
Before we define all the objects more precisely, it may help the reader
to have a loose description of the construction. We will define by
induction a sequence of quadruples
$(\bbf{e}^n,D_{n-1},\hat\cE_n,U^n)$ for $n\ge 1$, where:
\begin{itemize}
\item $\bbf{e}^n=(\bbf{a}^n,\bbf{b}^n)$ is a {\em frame}, that is,
a point of a certain closed subscheme $\Frame_{n-1}$ of a certain fibred product
$W^{n-1}\times (W^\lambda)^{n-1}$.
Frames are the parameters of the construction, to be chosen at each step.
\item $D_{n-1}:\hat\cE_{n-1}\to\hat\GG_m$ is a morphism of formal $R$-schemes
which mod $\lambda_n$ induces a morphism of formal $(R/\lambda_nR)$-groups.
\item $\hat\cE_n$ is a commutative formal group extension of
$\hat\cE_{n-1}$ by $\hat\cG^{\lambda_n}$ such that the map
$\alpha_{\hat\cE_n}:\hat\cE_n\to (\hat\GG_m)^n$ defined on the points by
$$
(x_1,\dots,x_n)\mapsto
(D_0+\lambda_1 x_1,D_1+\lambda_2 x_2,\dots,D_{n-1}+\lambda_n x_n)
$$
is a morphism of formal groups, where $D_i=D_i(x_1,\dots,x_i)$ for the
natural coordinates $x_1,\dots,x_i$ on $\hat\cE_i$.
\item $U^n:W^n\to W^n$ is a morphism of $R$-group schemes.
\end{itemize}
\end{noth}

\begin{noth} {\bf Initialization.}
The induction is initialized at $n=1$. Let $W^0=0$ and $\hat\cE_0=0$.
We set $\bbf{e}^1=(0,0)$, $D_0:\hat\cE_0\to\hat\GG_m$ equal
to $1$, $\hat\cE_1=\hat\cG^{\lambda_1}$ and $U^1=F^{\lambda_1}:W\to W$.
\end{noth}

\begin{noth} {\bf Induction.}
For the inductive step of the construction, we assume that
$(\bbf{e}^i,D_{i-1},\hat\cE_i,U^i)$ has been constructed for $1\le i\le n$ and
we explain how to produce $(\bbf{e}^{n+1},D_n,\hat\cE_{n+1},U^{n+1})$.
For this, we introduce
frames. Let $\lambda\in R$ be a nonzerodivisor and consider the morphism
$$U^n-\lambda:W^n\times(W^{\lambda})^n\to W^n$$
taking an element $(\bbf{a}^{n+1},\bbf{b}^{n+1})\in W^n\times(W^{\lambda})^n$
to $U^n(\bbf{a}^{n+1})-\lambda.\bbf{b}^{n+1}$.

\begin{defi} \label{df:frame}
A {\em $\lambda$-frame (relative to $\cE_n$)} is an $R$-point
$\bbf{e}^{n+1}=(\bbf{a}^{n+1},\bbf{b}^{n+1})$ of the kernel of
$U^n-\lambda$. The {\em scheme of frames of dimension $n$} is
$\Frame_n=\ker(U^n-\lambda)$.
\end{defi}

Now the induction goes in four steps A-B-C-D.

\bigskip

\no {\bf A.} Choose a $\lambda_{n+1}$-frame
$\bbf{e}^{n+1}=(\bbf{a}^{n+1},\bbf{b}^{n+1})\in \Frame_n(R)$.

\bigskip

\no {\bf B.} It is in the definition and properties of $D_n$ that
lies the main input of Sekiguchi-Suwa theory.
Let $A$ be an $R$-algebra. Let us extend the terminology of
Definition~\ref{df:fund_1} by calling a morphism of formal $A$-schemes
$\hat\cE_{n,A}\to\hat\GG_{m,A}$ {\em fundamental} if it is a product
of Artin-Hasse exponentials
$$
E_p(\bbf{a}^{n+1}_1,\lambda_1,X_1/D_0)\,E_p(\bbf{a}^{n+1}_2,\lambda_2,X_2/D_1)
\dots E_p(\bbf{a}^{n+1}_n,\lambda_n,X_n/D_{n-1})
$$
for some $n$-tuple of Witt vectors
$\bbf{a}^{n+1}=(\bbf{a}^{n+1}_1,\dots,\bbf{a}^{n+1}_n)\in W(A)^n$.
Then, we have:

\begin{theo} \label{th:AHE_n}
Denote by $\FS/R$ the category of formal $R$-schemes and by
$\FG/R$ the category of formal $R$-groups. Then with the above notation
we have:
\begin{trivlist}
\itemn{1}
The deformed Artin-Hasse exponentials define a
monomorphism of $R$-group functors
$$
\Fund:W^n\too\Hom_{\FS/R}(\hat\cE_n,\hat\GG_m)
$$
taking an $n$-tuple of Witt vectors
$\bbf{a}^{n+1}=(\bbf{a}^{n+1}_1,\dots,\bbf{a}^{n+1}_n)\in W(A)^n$
to the corresponding fundamental morphism
$\prod_{i=1}^n\,E_p(\bbf{a}^{n+1}_i,\lambda_i,X_i/D_{i-1})$.
Here, the group law on the target is induced by the group law of $\hat\GG_m$.
\itemn{2}
The map $\Fund$ induces an isomorphism of $R$-group functors
$$
\ker(U^n:W^n\to W^n)\stackrel{\sim}{\too}\Hom_{\FG/R}(\hat\cE_n,\hat\GG_m).
$$
In particular, any morphism of formal $R$-groups
$\hat\cE_n\to \hat\GG_m$ is fundamental.
\end{trivlist}
\end{theo}

\begin{proo}
Point (1) is \cite{SS1}, Corollary~2.9 and point (2) is \cite{SS2}, Theorem~5.1.
\end{proo}

It follows from the definition of a frame and from point (2) of the
theorem that if we take for $\bbf{a}^{n+1}$ the first component of the
frame $\bbf{e}^{n+1}=(\bbf{a}^{n+1},\bbf{b}^{n+1})$ chosen in Step~A,
then $\bbf{a}^{n+1}$ lies in the kernel of $U^n$ modulo $\lambda_{n+1}$
and the fundamental morphism of formal $R$-schemes
$$
D_n=\prod_{i=1}^n\,E_p(\bbf{a}^{n+1}_i,\lambda_i,X_i/D_{i-1})
$$
induces modulo $\lambda_{n+1}$ a morphism of formal $(R/\lambda_{n+1}R)$-groups.

\bigskip

\no {\bf C.} We now build $\hat\cE_{n+1}$.
Since $D_n$ gives a morphism of formal $(R/\lambda_{n+1}R)$-groups,
then the expression $D_n(X)D_n(Y)D_n(X\star Y)^{-1}-1$ vanishes mod $\lambda_{n+1}$,
where $X\star Y$ denotes the group law in $\hat\cE_n$. Since $\lambda_{n+1}$
is a nonzerodivisor, this implies that
$$
H_n(X,Y)=\frac{1}{\lambda_{n+1}}\left(\frac{D_n(X)D_n(Y)}{D_n(X\star Y)}-1\right)
$$
is well-defined. It is a symmetric $2$-cocycle
$\hat\cE_n\times\hat\cE_n\to\hat\cG^{\lambda}$ i.e. an element of the
Hochschild cohomology group
$\HH^2_0(\hat\cE_n,\hat\cG^{\lambda})$ of symmetric $2$-cocycles.
From a $2$-cocycle we can construct an extension of $\hat\cE_n$ by
$\hat\cG^{\lambda}$ in the usual way: this is $\hat\cE_{n+1}$.

\bigskip

\no {\bf D.} Define $U^{n+1}:W^{n+1}\to W^{n+1}$ by the matrix
$$
U^{n+1}=\left(
\begin{array}{cccc}
& & & -T_{\bbf{b}^{n+1}_1} \\
& U^n & & \vdots \\
& & & -T_{\bbf{b}^{n+1}_n} \\
0 & \dots & 0 & F^{\lambda_{n+1}} \\
\end{array}
\right).
$$

\end{noth}

With the following definition and theorem, we point out that this
construction is universal:

\begin{defi}
A {\em framed formal $R$-group of type $(\lambda_1,\dots,\lambda_n)$}
is a sequence
$$
\hat\cE_0=0,(\hat\cE_1,\bbf{e}^1),\dots,(\hat\cE_n,\bbf{e}^n)
$$
of pairs composed of an affine smooth commutative formal group scheme
and a frame, such that for each $i=1,\dots,n$ the  formal group scheme
$\hat\cE_i$ is the extension of $\hat\cE_{i-1}$ by $\hat\cG^{\lambda_i}$
determined by the $\lambda_i$-frame $\bbf{e}^i$.
We often write $\hat\cE_n$ as a shortcut for this data.
\end{defi}

\begin{theo} \label{th:universal_formal}
Let $\AA^n=\Spec(\ZZ_{(p)}[\Lambda_1,\dots,\Lambda_n])$ be affine
$n$-space over $\ZZ_{(p)}$. Then there exists an affine flat
$\AA^n$-scheme $\sS_n=\Spec(\sR_n)$ and a framed formal $\sR_n$-group
$\hat\sE_n$ of type $(\Lambda_1,\dots,\Lambda_n)$ with the
following universal property~: for any $\ZZ_{(p)}$-algebra $R$, any
nonzerodivisors $\lambda_1,\dots,\lambda_n\in R$ and any framed formal
$R$-group $\hat\cE_n$ of type $(\lambda_1,\dots,\lambda_n)$,
there exists a unique map $\sR_n\to R$ taking $\Lambda_i$ to $\lambda_i$
such that $\hat\cE_n\simeq \hat\sE_n\otimes_{\sR_n} R$.
\end{theo}

\begin{proo}
The proof is almost tautological, because framed formal groups are
more or less by construction pullback of a universal one. Let us however
sketch it. What we have to do is to carry out the induction as before,
in a universal way. Let $W^0=0$ and $\hat\cE_0=0$.

For $n=1$ we put $\sR_1=\ZZ_{(p)}[\Lambda_1]$, $\bbf{e}^1=(0,0)$,
$D_0=1$, $\hat\sE_1=\hat\cG^{\Lambda_1}$ and $U^1=F^{\Lambda_1}:W\to W$.

Once $\sS_i$, $\bbf{e}^i$, $D_{i-1}$, $\hat\sE_i$ and $U^i$ have been
constructed for $1\le i\le n$, we find $\sS_{n+1}$, $\bbf{e}^{n+1}$,
$D_n$, $\hat\sE_{n+1}$ and $U^{n+1}$ as follows. We take as a base ring
the ring $R':=\sR_n\otimes\ZZ_{(p)}[\Lambda_{n+1}]$. We define $\sS_{n+1}$
as the scheme of frames $\Frame_n=\ker(U^n-{\Lambda_{n+1}})$,
and we set $\bbf{e}^{n+1}$ equal to the universal point of $\sS^{n+1}$.
Note that since $U^n$ is given by a triangular matrix whose diagonal
entries are flat morphisms by Lemma~\ref{lm:F_flat}, it follows
immediately that it is a flat morphism. By the definition of $\sS_{n+1}$
as the fibred product
$$
\xymatrix{
\sS_{n+1} \ar[r] \ar[d] & (W^{\Lambda_{n+1}})^n \ar[d]^{{\Lambda_{n+1}}} \\
W^n \ar[r]^{U^n} & W^n \\}
$$
we see that it is flat over $(W^{\Lambda_{n+1}})^n$, hence flat over
$\AA^{n+1}$. It follows that $\Lambda_{n+1}$ is not a zerodivisor in the
function ring $\sR_{n+1}$ of $\sS_{n+1}$. Now the coefficient $\bbf{a}^{n+1}$
of the frame $\bbf{e}^{n+1}$ determines a fundamental morphism
$$D_n=\prod_{i=1}^n\,E_p(\bbf{a}^{n+1}_i,\Lambda_i,X_i/D_{i-1}),$$
a $2$-cocyle
$$
H_n(X,Y)=\frac{1}{\Lambda_{n+1}}\left(\frac{D_n(X)D_n(Y)}{D_n(X\star Y)}-1\right)
$$
and then an extension $\hat\sE_{n+1}$ in the same way as before. The
coefficient $\bbf{b}^{n+1}$ of the frame determines a matrix $U^{n+1}$
by the same formula as in Step~D of the induction. Once the construction
is over, the verification of the universal property is immediate.
\end{proo}

\section{Algebraic theory} \label{sec:AT}

In this section, we show how to adapt the formal constructions in order
to provide (algebraic) filtered group schemes. This is done by
truncating the power series and the Witt vector coefficient in a
suitable way. We give some preliminaries on truncations in
Subsections~\ref{truncation_dahes} and~\ref{truncation_WV}.
Then we proceed to construct filtered
group schemes in~\ref{ss:fil_gs}, with Theorem~\ref{th:universal_alg}
as the final point.

\subsection{Truncation of deformed Artin-Hasse exponentials}
\label{truncation_dahes}

In order to produce non-formal group schemes, we will need the deformed
exponentials to be polynomials. We can achieve this either by letting enough
coefficients specialize to nilpotent elements, or by truncating.
We know from \cite{SS2}, Prop. 2.11 that if $\Lambda,U_0,U_1,\dots$ specialize
to nilpotent elements, only finitely many of them nonzero, then
$E_p(\UU,\Lambda,T)$ specializes to a polynomial. In the following lemma,
we give an exact bound for the degree of this polynomial, in terms of bounds
on the number of nonzero coefficients and the nilpotency indices.

\begin{lemm} \label{lm:TDE}
Let $L,M,N\ge 1$ be integers. Then if we reduce the coefficients of the
deformed exponential $E_p(\UU,\Lambda,T)$ modulo the ideal generated by
$$
\Lambda^L, (U_0)^N,(U_1)^N,\dots,(U_{M-1})^N,U_M,U_{M+1},\dots
$$
then the series $E_p(\UU,\Lambda,T)$ is a polynomial of degree at most
$(N-1)\frac{p^M-1}{p-1}+(L-1)$.
\end{lemm}

\begin{proo}
For each $\ell$, we have
$E_p(U_\ell,\Lambda^{p^\ell},T^{p^\ell})=
E_p(U_\ell T^{p^\ell},\Lambda^{p^\ell} T^{p^\ell},1)$.
It follows that the latter series is a sum of monomials of the form
$(U_\ell T^{p^\ell})^i(\Lambda^{p^\ell}T^{p^\ell})^j$ for varying $i,j$.
Now let us take images in the indicated quotient ring. There, for all
$\ell\ge M$ we have $U_\ell=0$ and $E_p(U_\ell,\Lambda^{p^\ell},T^{p^\ell})=1$.
It follows that only the first $M$ factors show up in the product defining
$E_p(\UU,\Lambda,T)$. A typical monomial in this series is obtained
by picking a monomial of index $i_\ell,j_\ell$ in each factor; the result is
the product of
$$
(U_0)^{i_0}(U_1)^{i_1}\dots (U_{M-1})^{i_{M-1}}
\times
T^{i_0+i_1p+\dots+i_{M-1}p^{M-1}}
$$
by
$$
\Lambda^{j_0+j_1p+\dots+j_{M-1}p^{M-1}}
\times
T^{j_0+j_1p+\dots+j_{M-1}p^{M-1}}.
$$
For this to be nonzero, we must have $i_\ell\le N-1$ for each $\ell$ and
$$j_0+j_1p+\dots+j_{M-1}p^{M-1}\le L-1$$
for each $(j_0,\dots,j_{M-1})$.
Thereby the $T$-degree of the monomial is less than
$$
(N-1)(1+p+\dots+p^{M-1})+(L-1),
$$
which is what the lemma claims.
\end{proo}

\begin{defi} \label{df:TDE}
Let $L,M,N\ge 1$ be integers and let $\tau_{L,M,N}$ be the truncation map of
power series in degrees $\ge (N-1)\frac{p^M-1}{p-1}+(L-1)+1$. Then the polynomial
$$
E^{L,M,N}_p(\UU,\Lambda,T)\egaldef\tau_{L,M,N} E_p(\UU,\Lambda,T)
\in \ZZ_{(p)}[\Lambda,U_0,U_1,\dots][T]
$$
is called the {\em truncated (deformed) exponential of level $(L,M,N)$}.
\end{defi}

\subsection{Truncation of Witt vectors}
\label{truncation_WV}

We will make big use of the functor $\hat{W}$ and its pushforward
$i_*\hat W$ by the closed immersion $i:\Spec(\ZZ)\into \AA^1=\Spec(\ZZ[\Lambda])$.
Since $\hat W$ is naturally filtered, this leads to consider various
truncations of $W$ and $\hat{W}$, over $\Spec(\ZZ)$ and over $\AA^1$.
In order to define them, we fix integers $M,N\ge 1$.

\begin{noth}
{\bf Truncation by the length.}
\begin{trivlist}
\itemn{1} $W_M$ is the $\ZZ$-subfunctor of $W$ defined by
$W_M(A)=\{a\in W(A),\,a_i=0 \mbox{ for } i\ge M\}$. We emphasize that it is
of course {\em not} a subgroup functor; it should not be confused with the
quotient ring of Witt vectors of length $N$, which will not appear in the
present paper.
\itemn{2} $W^\Lambda_M$ is the $\AA^1$-subfunctor of $W^\Lambda$ defined by
$W^\Lambda_M(A)=\{a\in W^\Lambda(A),\,a_i=0 \mbox{ for } i\ge M\}$.
\itemn{3} $\hat W_M=\hat W\cap W_M$ is a $\ZZ$-subfunctor of $\hat W$.
\end{trivlist}
\end{noth}

\begin{noth}
{\bf Truncation by the nilpotency index.}
\begin{trivlist}
\itemn{4} $W_{M,N,\Lambda}\subset W_M$ is the $\AA^1$-subfunctor
defined by
$$
W_{M,N,\Lambda}(A)=\{\bbf{a}\in W_M(A),\, (a_i)^N\equiv 0 \mod \Lambda
\mbox{ for all } i\}.
$$
\itemn{5} $\hat W_{M,N}=W_{M,N,0}\subset \hat W_M$ is the
$\ZZ$-subfunctor of $\hat W_M$
introduced in the proof of Lemma~\ref{lm:ideal}.
\end{trivlist}
\end{noth}

We view all these functors as sheaves over the small flat sites
$\Spec(\ZZ)_{\ffl}$ and $\AA^1_{\ffl}$. Then $W_M$ and $W^\Lambda_M$
are representable by $M$-dimensional affine spaces over $\Spec(\ZZ)$,
$\hat W_{M,N}$ is representable by a finite flat $\ZZ$-scheme, and
$W_{M,N,\Lambda}$ is representable by a scheme which is a finite flat
$N^M$-sheeted cover of an $M$-dimensional affine space over $\AA^1$.
Of these statements, only the last deserves a comment. The basic
observation is that the sheaf $F$ on $\AA^1_{\ffl}$ defined by
$F(A)=\{a\in A,\,a^N\equiv 0 \mod\Lambda\}$ is represented by the scheme
$\Spec(\ZZ[\Lambda][u,v]/(u^N-\Lambda v))$, and then $W_{M,N,\Lambda}$
is obviously represented by the $M$-fold product of $F$.

\subsection{Construction of framed group schemes} \label{ss:fil_gs}

Let $R$ be a $\ZZ_{(p)}$-algebra and $\lambda_1,\lambda_2,\dots$ elements of $R$.
Filtered $R$-group schemes are defined just like their formal analogues
in Definition~\ref{df:FFGS}.

\begin{defi} \label{df:FiltGS}
A {\em filtered $R$-group scheme of type $(\lambda_1,\dots,\lambda_n)$}
is a sequence
$$\cE_0=0,\cE_1,\dots,\cE_n$$
of affine smooth commutative group schemes such that for each
$i=1,\dots,n$ the group scheme $\cE_i$ is an extension of
$\cE_{i-1}$ by $\cG^{\lambda_i}$.
\end{defi}

\begin{assu} \label{ass:radicals}
The elements $\lambda_1,\lambda_2,\dots$ are not zero divisors in $R$, and
$\lambda_i$ is nilpotent modulo $\lambda_{i+1}$ for each $i\ge 1$.
\end{assu}

We will see that under this assumption,
and provided we make suitable truncations, the procedure described in
Subsection~\ref{ss:fil_fgs} in the formal case gives filtered group schemes.
In order to carry out the construction, we fix positive integers
$L_1,L_2,\dots$ such that $(\lambda_i)^{L_i}\in \lambda_{i+1}R$ for all $i\ge 1$.
We also fix a pair of positive integers $(M,N)$ serving as a truncation level.

\begin{noth} {\bf Description of the procedure.} \label{no:desc_proc}
Contrary to the formal situation, here the polynomials giving the fundamental
morphisms $D_i$ will not be invertible over $R$ but only over
$R/\lambda_{i+1}$. Because the inductive definition of the $D_i$ requires
lifts of the inverses, we have to consider such lifts to be part of
the data that we need to produce. So this time, the $n$-th step of the
induction will produce $5$-tuples
$(\bbf{e}^n,D_{n-1},D_{n-1}^{-1},\cE_n,U^n)$ where, more precisely:
\begin{itemize}
\item $\bbf{e}^n=(\bbf{a}^n,\bbf{b}^n)$ is a {\em frame}, that is,
a point of a certain closed subscheme $\Frame_{n-1}$ of a certain fibred product
of $(W_{M,N,\lambda_n})^{n-1}$ by $(W_{M_n}^{\lambda_n})^{n-1}$.
Frames are the parameters of the construction, to be chosen at each step.
\item $D_{n-1},D_{n-1}^{-1}:\cE_{n-1}\to\AA^1$ are truncated deformed
exponentials, that is morphisms of $R$-schemes
which mod $\lambda_n$ induce mutually inverse morphisms of
$(R/\lambda_nR)$-group schemes $\cE_{n-1}\to\GG_m$.
\item $\cE_n$ is a commutative $R$-group scheme extension of
$\cE_{n-1}$ by $\cG^{\lambda_n}$, with underlying scheme
$$
\cE_n=\Spec\left( R\left[X_1,\dots,X_n,\frac{1}{D_0+\lambda_1 X_1},\dots,
\frac{1}{D_{n-1}+\lambda_n X_n}\right]\right),
$$
such that the map $\alpha_{\cE_n}:\cE_n\to (\GG_m)^n$ defined on the points by
$$
(x_1,\dots,x_n)\mapsto
(D_0+\lambda_1 x_1,D_1+\lambda_2 x_2,\dots,D_{n-1}+\lambda_n x_n)
$$
is a morphism of $R$-group schemes.
\item $U^n:(W_{M,N,\lambda_{n+1}})^n\to (W_{M_n,N_n,\lambda_{n+1}})^n$
is a morphism of $R$-schemes represented by a square matrix of size
$n$, where $M_n,N_n$ are integers.
\end{itemize}
\end{noth}

\begin{noth} {\bf Initialization.} \label{init}
We set $W^0=(W_{M,N})^0=0$ and $\cE_0=0$.
The induction is initialized at $n=1$ by setting $\bbf{e}^1=(0,0)$,
$D_0=D_0^{-1}:\cE_0\to\GG_m\subset\AA^1$ equal to $1$, and
$\cE_1=\cG^{\lambda_1}$. It follows from Lemmas~\ref{lm:ideal} and
\ref{lm:F_stable} that the endomorphism
$F^{\lambda_1}:W\otimes (R/\lambda_2)\to W\otimes (R/\lambda_2)$
leaves $\hat W\otimes (R/\lambda_2)$ stable, so it maps
$\hat W_{M,N}\otimes (R/\lambda_2)$ into
$\hat W_{M_1,N_1}\otimes (R/\lambda_2)$ for some integers
$M_1$, $N_1$. It follows that the composition of
$F^{\lambda_1}:W_{M,N,\lambda_2}\to W$ with the truncation
map $\tau_{\ge M_1}:W\to W_{M_1}$ factors through
$W_{M_1,N_1,\lambda_2}$. The result is a morphism
$U^1=F^{\lambda_1}:W_{M,N,\lambda_2}\to W_{M_1,N_1,\lambda_2}$.
\end{noth}

\begin{noth} {\bf Induction.} \label{induct}
For the inductive step of the construction, we assume that
$$(\bbf{e}^i,D_{i-1},D_{i-1}^{-1},\cE_i,U^i)$$ has been constructed
for $1\le i\le n$ and we explain how to produce
$(\bbf{e}^{n+1},D_n,D_n^{-1},\cE_{n+1},U^{n+1})$. We do this in four
steps A-B-C-D.

\medskip

\no {\bf A.}
To start with, we choose $\bbf{e}^{n+1}=(\bbf{a}^{n+1},\bbf{b}^{n+1})$
such that $U^n(\bbf{a}^{n+1})=\lambda_{n+1}.\bbf{b}^{n+1}$.
To be more formal, this is a section over $R$ of the {\em scheme of frames}
$\Frame_n$ defined as the fibred product of the morphisms
$$
U^n:(W_{M,N,\lambda_{n+1}})^n\to
(W_{M_n,N_n,\lambda_{n+1}})^n\subset(W_{M_n})^n
\quad \mbox{ and } \quad
{\lambda_{n+1}}:(W_{M_n}^{\lambda_{n+1}})^n\to (W_{M_n})^n,
$$
that is:
$$
\Frame_n=(W_{M,N,\lambda_{n+1}})^n\times_{(W_{M_n})^n}
(W_{M_n}^{\lambda_{n+1}})^n.
$$
The choice of $\bbf{e}^{n+1}$ will determine the other four objects in
the $5$-tuple.

\medskip

\no {\bf B.}
Using the first component $\bbf{a}^{n+1}$ of the frame, we define:
$$
\begin{array}{lcl}
D_n & =
& \prod_{i=1}^n\,E^{L_i,M_i,N_i}_p(\bbf{a}^{n+1}_i,\lambda_i,D_{i-1}^{-1}X_i)
\bigskip \\
D_n^{-1}& =
& \prod_{i=1}^n\,E^{L_i,M_i,N_i}_p(-\bbf{a}^{n+1}_i,\lambda_i,D_{i-1}^{-1}X_i). \\
\end{array}
$$
Note that this is where the $D_i^{-1}$ are useful, since they are involved in
the definition of the $D_i$.
It follows from the choice of the truncations (involved in the choice
of $\bbf{a}^{n+1}$ and in the truncated exponentials, see Lemma~\ref{lm:TDE}
and Definition~\ref{df:TDE})
and from Theorem~5.1 of \cite{SS2}
(in the case of nilpotent coefficients), that $D_n$ and $D_n^{-1}$
induce morphisms of $(R/\lambda_{n+1}R)$-group schemes $\cE_n\to \GG_m$
inverse to each other.

\medskip

\no {\bf C.}
Now we define $\cE_{n+1}$. At this step, the strategy differs from the
formal case because the truncated deformed exponentials are not invertible
and do not give rise to $2$-cocycles like in the formal case. In fact,
the Hochschild cohomology group $\HH^2_0(\cE_n,\cG^{\lambda_{n+1}})$
is usually very small. Instead, we use the exact sequence of sheaves
on the small flat site
$$
0\too\cG^{\lambda_{n+1}}\too\GG_m\too i_*\GG_m\too 0
$$
where $i:\Spec(R/\lambda_{n+1}R)\into\Spec(R)$ is the closed immersion.
There is a coboundary map
$\Hom(\cE_n,i_*\GG_m)\to\Ext^1(\cE_n,\cG^{\lambda_{n+1}})$ and the extension
$\cE_{n+1}$ is the image of $D_n:\cE_n\to i_*\GG_m$ under this map.
(Note that all the cohomology groups may be computed in the big flat site,
by Milne \cite{Mi}, III, Remark~3.2.) It is
the extension obtained by pulling back the extension
$0\to\cG^{\lambda_{n+1}}\to\GG_m\to i_*\GG_m\to 0$
along $D_n$. Thus for each flat $R$-algebra $A$, we have:
$$
\begin{array}{rl}
\cE_{n+1}(A)
& =\{(v,w)\in \cE_n(A)\times A^\times\,,\,D_n(v)\equiv w
\mod \lambda_{n+1}\} \medskip \\
& =\{(v,w)\in \cE_n(A)\times A^\times\,,\,D_n(v)+\lambda_{n+1} x=w
\mbox{ for some } x\in A\} \medskip \\
& =\{(v,x)\in \cE_n(A)\times A\,,\,D_n(v)+\lambda_{n+1} x\in A^\times\}.
\end{array}
$$
This sheaf is represented by the scheme
$\cE_{n+1}=\Spec(R[\cE_n][X_n,\frac{1}{D_n+\lambda_{n+1}
X_{n+1}}])$. As far as the group law is concerned, note that by
the assumption on $D_n$ there exists a unique function $K=K(X,Y)$
on $\cE_n\times\cE_n$ such that $D_n(X)D_n(Y)=D_n(X\star
Y)+\lambda_{n+1}K(X,Y)$, where $X\star Y$ denotes the group law in
$\cE_n$. Then it is easy to see that the group law in $\cE_{n+1}$
is given on the points by:
$$
(v_1,x_1)\star' (v_2,x_2)
=\big(v_1\star v_2,x_1D_n(v_2)+x_2D_n(v_1)+\lambda_{n+1}x_1x_2+K(v_1,v_2)\big).
$$
Equivalently, the group law is the only one such that the map
$$
\alpha_{\cE_{n+1}}:
\begin{array}[t]{l}
\cE_{n+1}\to (\GG_m)^{n+1} \\
(x_1,\dots,x_{n+1})\mapsto
(D_0+\lambda_1 x_1,D_1+\lambda_2 x_2,\dots,D_n+\lambda_{n+1} x_{n+1}) \\
\end{array}
$$
is a morphism of $R$-group schemes.

\medskip

\no {\bf D.}
Finally, using the second component $\bbf{b}^{n+1}$ of the frame,
we consider the matrix:
$$
U^{n+1}=\left(
\begin{array}{cccc}
& & & -T_{\bbf{b}^{n+1}_1} \\
& U^n & & \vdots \\
& & & -T_{\bbf{b}^{n+1}_n} \\
0 & \dots & 0 & F^{\lambda_{n+1}} \\
\end{array}
\right).
$$
Let $\lambda\in R$ be a nonzerodivisor such that $\lambda_{n+1}$ is nilpotent
modulo $\lambda$. If we reduce modulo $\lambda$, then according to
Lemmas~\ref{lm:ideal}, \ref{lm:F_stable}, \ref{lm:T_stable},
the endomorphism
$$
U^{n+1}\otimes (R/\lambda R):W^{n+1}\otimes (R/\lambda R)
\to W^{n+1}\otimes (R/\lambda R)
$$
leaves $\hat W^{n+1}\otimes (R/\lambda R)$ stable. It follows that
there exist integers $M_{n+1}$, $N_{n+1}$ such that
$U^{n+1}\otimes (R/\lambda R)$ maps $(\hat W_{M,N})^{n+1}$
into $(\hat W_{M_{n+1},N_{n+1}})^{n+1}$. Therefore the composition of
$$
U^{n+1}:(W_{M,N,\lambda})^{n+1}\subset W^{n+1} \to W^{n+1}
$$
with the truncation map $W^{n+1}\to (W_{M_{n+1}})^{n+1}$
factors through $(W_{M_{n+1},N_{n+1},\lambda})^{n+1}$.
Fixing $\lambda=\lambda_{n+2}$, the result is a morphism
of $R$-schemes
$$
U^{n+1}:(W_{M,N,\lambda_{n+2}})^{n+1}
\to (W_{M_{n+1},N_{n+1},\lambda_{n+2}})^{n+1}
$$
This is the last object in our sought-for $5$-tuple.
\end{noth}

\begin{rema}
A priori, the integers $M_n,N_n$ depend on the particular frames involved
in the matrices $U^n$. However,
considering the universal case (see Theorem~\ref{th:universal_alg}), it is
seen that in fact, once $(M,N)$ is fixed then $(M_n,N_n)$ may be chosen
uniform, minimal and hence completely determined by $M_1,N_1$ and $n$.
\end{rema}

\begin{defi}
A {\em framed $R$-group scheme of type $(\lambda_1,\dots,\lambda_n)$}
is a sequence
$$
\cE_0=0,(\cE_1,\bbf{e}^1),\dots,(\cE_n,\bbf{e}^n)
$$
of pairs composed of an affine smooth commutative group scheme
and a frame, such that for each $i=1,\dots,n$ the group scheme
$\cE_i$ is the extension of $\cE_{i-1}$ by $\cG^{\lambda_i}$
determined by the frame~$\bbf{e}^i$.
We often write $\cE_n$ as a shortcut for this data.
\end{defi}

In order to state the analogue of Theorem~\ref{th:universal_formal} in
the algebraic context, we must make sure that the coefficients $\lambda_i$
satisfy Assumption~\ref{ass:radicals}. This means that for some integer
$\nu\ge 1$ they are points of the space
$$
\BB^n_\nu=\Spec\left(\frac{\ZZ_{(p)}[\Lambda_1,\dots,\Lambda_n,M_2,\dots,M_n]}
{\Lambda_1^\nu-M_2\Lambda_2,\dots,\Lambda_{n-1}^\nu-M_n\Lambda_n}\right) \, .
$$
This is a finite flat cover of the affine space
$\AA^n=\Spec(\ZZ_{(p)}[M_2,\dots,M_n,\Lambda_n])$. Moreover, there are
obvious projections $\BB^{n+1}_\nu\to\BB^n_\nu$ given by the inclusion
of function rings. Below, we denote by $\Lambda$ the product of the $\Lambda_i$.

\begin{theo} \label{th:universal_alg}
Let $\BB^n_\nu$ be the finite flat covers of affine space $\AA^n$
defined above. Then there exists a sequence indexed by $\nu\ge 1$ of affine
$\BB^n_\nu$-schemes $\sS_n^\nu=\Spec(\sR_n^\nu)$ of finite type,
without $\Lambda$-torsion, and framed
$\sR_n^\nu$-group schemes $\sE_n^\nu$ of type $(\Lambda_1,\dots,\Lambda_n)$
with the following universal property~: for any $\ZZ_{(p)}$-algebra $R$,
any nonzerodivisors $\lambda_1,\dots,\lambda_n\in R$ such that $\lambda_i$
is nilpotent modulo $\lambda_{i+1}$ for each $i$, and any framed
$R$-group scheme $\cE_n$ of type $(\lambda_1,\dots,\lambda_n)$, there
exists $\nu$ and a unique map $\sR_n^\nu\to R$ taking $\Lambda_i$ to
$\lambda_i$ such that $\cE_n\simeq \sE_n^\nu\otimes_{\sR_n^\nu} R$.
\end{theo}

\begin{proo}
For a fixed $\nu \ge 1$, we first give $\sS_n^\nu\to \BB^n_\nu$ and
$\sE_n^\nu\to \sS_n^\nu$. The construction is by induction on $n$
and follows the proof of Theorem~\ref{th:universal_formal}, with
minor differences which we indicate. The main difference is that in
the present case, the function ring of the schemes of frames in
dimension $n$ is bound to play the role of the coefficient ring in
dimensions $\ge n+1$ and so needs to be free of $\Lambda$-torsion.
Thus we have to kill torsion in the adequate fibred product.

We set $L=M=N=\nu$. In order to keep the notation light, we will sometimes
omit the symbol $\nu$ in the indices and exponents.
We initialize by $\sE_0=0$, $\sS_1=\BB^1$,
$\bbf{e}^1=(0,0)$, $D_0=D_0^{-1}=1$,
$\sE_1=\cG^{\Lambda_1}$ over $\sS_1$, and
$U^1=F^{\Lambda_1}:W_{M,N,\Lambda_2}\to W_{M_1,N_1,\Lambda_2}$
is the morphism of $\BB^2$-schemes constructed  like in~\ref{init}.
Now assuming that for $1\le i\le n$ we have objects
$\sS_i,\bbf{e}^i,D_{i-1},D_{i-1}^{-1},\sE_i,U^i$, here is how to
construct $\sS_{n+1},\bbf{e}^{n+1},D_n,D_n^{-1},\sE_{n+1},U^{n+1}$.

Consider the morphisms of $\sS_n\times_{\BB^n}\BB^{n+1}$-schemes
$$
U^n:(W_{\nu,\nu,\Lambda_{n+1}})^n
\to (W_{M_{n+1},N_{n+1},\Lambda_{n+1}})^n\subset(W_{M_n})^n
\quad\mbox{ and } \quad
{\Lambda_{n+1}}:(W_{M_n}^{\Lambda_{n+1}})^n\to (W_{M_n})^n.
$$
Call $\sS_{n+1}$ the closed subscheme of the fibred product of $U^n$
and ${\Lambda_{n+1}}$ defined by the ideal of $\Lambda$-torsion,
where $\Lambda=\Lambda_1\dots\Lambda_{n+1}$. Let
$\bbf{e}^{n+1}=(\bbf{a}^{n+1},\bbf{b}^{n+1})$ be the universal point
of $\sS_{n+1}$. Then $D_n,D_n^{-1},\sE_{n+1},U^{n+1}$ are constructed
as in steps B, C, D of \ref{induct} and we do not repeat the details.

If $\cE_n$ is a framed group scheme of type $(\lambda_1,\dots,\lambda_n)$
over a ring $R$, then there exists $L$ such that
$(\lambda_i)^L \in \lambda_{i+1}R$. Moreover $\cE_n$ is described by
Witt vectors with a number of nonzero coefficients bounded by some $M$
and nilpotency indices bounded by some $N$. For $\nu=\max(L,M,N)$ it
is clear that $\cE_n$ is uniquely a pullback of $\sE_n^\nu$. This proves the
universality property of the statement of the theorem.
\end{proo}

\begin{prop}
Let $R$ be a $\ZZ_{(p)}$-algebra which is a unique factorization domain.
Then, any filtered group scheme is induced by a framed group scheme.
\end{prop}

\begin{proo}
By induction, it is enough to prove that given a filtered group scheme
$\cE$ of some type $(\lambda_1,\dots,\lambda_n)$ and a nonzero element
$\lambda\in R$, any extension of $\cE$ by $\cG^\lambda$ may be defined
by a frame. Consider the long exact sequence
$$
\dots\too\Hom(\cE_{R/\lambda},\GG_{m,R/\lambda})
\stackrel{\partial}{\too}\Ext^1(\cE,\cG^{\lambda})
\too\Ext^1(\cE,\GG_m)\too\dots
$$
derived from the exact sequence (1) in Proposition~\ref{pp:cG_cW}. It
is enough to prove that the coboundary $\partial$ is surjective. But this
follows from the fact that $\Ext^1(\cE,\GG_m)=0$, proven as in
\cite{SS2}, Proposition~3.1.
\end{proo}

\section{Kummer subschemes} \label{sec:KS}

Let $R$ be a $\ZZ_{(p)}$-algebra and let
$(\lambda_1,\dots,\lambda_n)$ be as in
Assumption~\ref{ass:radicals}.
We call $\lambda$ the product of the $\lambda_i$ and we write $K=R[1/\lambda]$.
For an $R$-scheme $X$, it will be a convenient abuse of terminology
to call the restriction $X_K$ the {\em generic fibre} of $X$.
Let $\cE$ be a framed group scheme of type $(\lambda_1,\dots,\lambda_n)$.
By construction $\cE$ comes with a map $\alpha_{\cE}:\cE\to (\GG_m)^n$ which
is an isomorphism over $K$. Let $\Theta^n:(\GG_m)^n\to (\GG_m)^n$ be the
morphism defined by
$$\Theta^n(t_1,\dots,t_n)=(t_1^p,t_2^pt_1^{-1},\dots,t_n^pt_{n-1}^{-1}).$$
The kernel of $\Theta^n$ is a subgroup isomorphic to $\mu_{p^n,R}$ which we
call the {\em Kummer $\mu_{p^n}$ of $\GG_m^n$}. Via the map $\alpha$, we
can see the Kummer $\mu_{p^n,K}$ as a closed subgroup scheme of
$\cE_K$. We define the {\em Kummer subscheme} as the scheme-theoretic
closure of $\mu_{p^n,K}$ in $\cE$. Note that in general the multiplication
of $G_K$ need not extend to $G$. The main question we want to address
in this section is~: when is the Kummer subscheme $G$ finite locally free
over $\Spec(R)$~? When this happens, then the multiplication extends and
accordingly, we shall prefer to call $G$ the {\em Kummer subgroup}.
In order to study this question, we first study the one-dimensional case
in~\ref{dim1}. Then, we consider extensions and we sketch the usual
inductive procedure producing isogenies between filtered group schemes,
in~\ref{dimn}.

Before we start, let us make a couple of easy remarks.
First, note that $G$ is the smallest
closed subscheme of $\cE$ with generic fibre isomorphic to
$\mu_{p^n,K}$.
It is also the only closed subscheme of $\cE$ without $\lambda$-torsion
with generic fibre isomorphic to $\mu_{p^n,K}$. In particular, if there
exists a closed subscheme of $\cE$ which is finite locally free over $R$
and has generic fibre isomorphic to $\mu_{p^n,K}$, then this subscheme
is equal to $G$.

\subsection{Dimension $1$} \label{dim1}

If $\lambda^{p-1}$ divides $p$ in $R$, then the polynomial
$\lambda^{-p}((\lambda x+1)^p-1)$ has coefficients in $R$ and
the morphism $\psi:\cG^\lambda\to\cG^{\lambda^p}$ defined
by $\psi(x)=\lambda^{-p}((\lambda x+1)^p-1)$ is an isogeny.
Following the notation in \cite{To}, we put $G_{\lambda,1}=\ker(\psi)$.

\begin{lemm} \label{lm:one_dim}
Let $\lambda\in R$ be a nonzerodivisor and $\cE=\cG^\lambda$.
\begin{trivlist}
\itemn{1} The Kummer subscheme $G$ is finite locally free over $R$
if and only if $\lambda^{p-1}$ divides $p$ in $R$.
\itemn{2} If $G$ is finite locally free, its ideal sheaf in $\cO_{\cE}$
is generated by
the polynomial $\lambda^{-p}((\lambda x+1)^p-1)$, and the quotient
$\cE\to \cE/G$ is isomorphic to the isogeny
$\psi:\cG^\lambda\to\cG^{\lambda^p}$.
\end{trivlist}
\end{lemm}

\begin{proo}
We begin with a couple of remarks.
Let us introduce the polynomial $P=(\lambda x+1)^p-1$. If
$s:=\max\,\{t\le p,\,\lambda^{t-1} \mbox{ divides } p\}$, there
exists $u\in R$ such that $p=u\lambda^{s-1}$.
Then we can write $P=\lambda^s Q$ where:
$$
Q=\lambda^{p-s}x^p+\sum_{i=1}^{p-1} \bin{p}{i} u \lambda^{i-1}x^i
$$
with ${p \choose i}=\bin{p}{i} p$ for $1\le i\le p-1$; and
$Q$ is not divisible by $\lambda$. The ideal of $G$ is:
$$
I=\big\{F\in R[x,(\lambda x+1)^{-1}],\,\exists n\ge 0,\,
\exists F'\in R[x,(\lambda x+1)^{-1}],\lambda^n F=QF'\big\}.
$$
Note that because $\lambda x+1$ is invertible modulo $P$ and also modulo $Q$,
we may always choose $F'\in R[x]$ above. Now let
$R[x]\to R[x,(\lambda x+1)^{-1}]$ be the natural inclusion and let
$J$ be the preimage of $I$. We have
$J=\{F\in R[x],\,\exists n\ge 0,\,\exists F'\in R[x],\lambda^n F=QF'\}$
and it is clear that the natural map $R[x]/J\to R[x,(\lambda x+1)^{-1}]/I$
is an isomorphism. We now prove (1) and (2).

\medskip

\no (1) We have to prove that the algebra $R[x]/J$ is finite
locally free over $R$ if and only if $\lambda^{p-1}$ divides $p$.
If $\lambda^{p-1}$ divides $p$, that is if $s=p$, then $Q$ is monic
and we claim that $J=(Q)$. Consider $F\in J$ and $n,F'$ such that
$\lambda^n F=QF'$.
We assume $n$ is minimal, i.e. $\lambda$ does not divide $F'$. If $n>0$
then $QF'\equiv 0 \mod \lambda$ hence $F'\equiv 0 \mod \lambda$ since $Q$
is monic hence a nonzerodivisor. This is a contradiction, so $n=0$
and $F\in (Q)$. Thus $J=(Q)$ and $R[x]/J$ is finite free over $R$.

Conversely, assume that $R[x]/J$ is finite locally free. We will prove
that $Q$ is monic and generates~$J$. It is enough to prove these properties
locally over $\Spec(R)$, hence we may assume that $R[x]/J$ is finite free
over $R$. Then there is a monic polynomial $G$ that generates $J$, see
Eisenbud \cite{Ei}, Prop.~4.1. From the fact that $Q\in (G)$ and
$\lambda^nG\in (Q)$
we see that $\deg(G)=\deg(Q)=p$. Writing $\lambda^nG=QF'$ we see that
$F'=\lambda^{n-p+s}$ so that $Q=\lambda^{p-s}G$.
Since $\lambda$ does not divide $Q$ this is possible only if $s=p$,
that is $\lambda^{p-1}$ divides $p$.

\medskip

\no (2) The isogeny $\psi:\cE=\cG^\lambda\to\cG^{\lambda^p}$
is $G$-invariant and induces a morphism $\cE/G\to\cG^{\lambda^p}$
which is finite flat of degree 1, hence an isomorphism.
\end{proo}

If $\cE$ is an $n$-dimensional framed group scheme, then what we have
just proved gives some "one-dimensional" necessary conditions for the
Kummer subscheme $G$ to be finite locally free over~$R$, as we shall
now see. Indeed if $G$ is finite locally free over $R$, then the quotient
$\cF=\cE/G$ is a smooth affine $n$-dimensional $R$-group scheme and the
quotient map $\nu:\cE\to\cF$ is an isogeny (smoothness follows from
\cite{EGA}, Chap.~0 (pr\'eliminaires), 17.3.3.(i)). Consider the subgroup
$\cG^{\lambda_n}\subset \cE$, its scheme-theoretic image $\cG$ under
$\nu$ and the restriction $\nu':\cG^{\lambda_n}\to\cG$ of $\nu$.
In the fibre over some point $s\in\Spec(R)$, the scheme $\cG_s$
is the quotient of $\cE_s$ by the equivalence relation induced by $G_s$,
that is, it is the quotient of $\cE_s$ by the stabilizer
$$
H=\{g\in G_s,\,g(\cG^{\lambda_n}_s)\subset \cG^{\lambda_n}_s\}.
$$
In particular $\cG_s$ is a quotient of a smooth $k(s)$-group scheme by a
finite flat subgroup scheme, hence it is a smooth $k(s)$-group scheme and
the map $\cG^{\lambda_n}_s\to \cG_s$ is flat. By the criterion of flatness
in fibres, it follows that $\nu'$ is flat and that $\cG$ is smooth.
Then the kernel $H_n=\ker(\nu')$ is flat of degree~$p$, with generic
fibre equal to the Kummer $\mu_{p,K}$ inside $G$. Moreover $\cG$ is
isomorphic to $\cG^{\lambda_n^p}$ and $\nu'$ is isomorphic to the isogeny
$\cG^{\lambda_n}\to\cG^{\lambda_n^p}$, by Lemma~\ref{lm:one_dim}.
Set $G_{n-1}=G/H_n$ and $\cF_{n-1}=\cF/\cG^{\lambda_n^p}$. Then we have
exact sequences
$$
0\too H_n\too G\too G_{n-1}\too 0,
$$
and
$$
0\too \cG^{p\lambda_n}\too\cF\too \cF_{n-1}\too 0.
$$
By induction we see immediately that $G$ and $\cF$ have filtrations
$G_0=0,G_1,\dots,G_n=G$ and $\cF_0=0,\cF_1,\dots,\cF_n=\cF$ where
$G_i\subset \cF_i$ is finite locally free of rank $p^i$ and
$\cF_i/\cF_{i-1}\simeq \cG^{\lambda_i^p}$. In particular $\cF$ is a
filtered group scheme of type $(\lambda_1^p,\dots,\lambda_n^p)$ and
$G$ is a successive extensions of the groups
$G_{\lambda_1,1},\dots,G_{\lambda_n,1}$. Another consequence of our
discussion is that the scheme-theoretic closure of $\mu_{p,K}$ inside
$\cG^{\lambda_n}$ is $H_n$ and in particular is finite locally free
over $R$. Similarly, by induction the scheme-theoretic closure of
$\mu_{p,K}$ inside $\cG^{\lambda_i}$ is equal to the kernel of
$G_i\to G_{i-1}$ and is finite locally free. By Lemma~\ref{lm:one_dim},
this proves that the following reinforcement of Assumption~\ref{ass:radicals}
is satisfied.

\begin{assu} \label{as:reinf}
For each $i\ge 1$ we have: $\lambda_i$ is not a zero divisor in $R$,
$\lambda_i$ is nilpotent modulo $\lambda_{i+1}$, and $\lambda_i^{p-1}$
divides $p$.
\end{assu}

\subsection{Construction of Kummer group schemes} \label{dimn}

From now on, we work under Assumption~\ref{as:reinf}. Because filtered
group schemes are defined by successive extensions, the condition that
the Kummer subscheme be finite locally free is also naturally expressed
at each extension step. Assume that we have a filtered group scheme $\cE_n$
of dimension $n$ with finite locally free Kummer subgroup $G_n$. Then, there
is a quotient morphism $\Psi_n:\cE_n\to\cF_n=\cE_n/G_n$ and for each
$\lambda\in R$ a pullback
$$
(\Psi^n)^*:\Hom_{R/\lambda R-\Gr}(\cF_n,\GG_m)\to
\Hom_{R/\lambda R-\Gr}(\cE_n,\GG_m).
$$
If $\cE_{n+1}$ is an extension of $\cE_n$ by $\cG^{\lambda}$
determined by a frame $\bbf{e}^{n+1}$, then we shall see that the condition
for the Kummer subscheme $G_{n+1}$ to be finite locally free is expressed
in terms of $(\Psi^n)^*$ and $\bbf{e}^{n+1}$. This will be integrated
in an inductive construction where we build at the same time the
group schemes $\cE_n$, $\cF_n$ and the isogeny between them, by making
compatible choices of frames. We explain how to do this, along the same
lines as before but giving a little less details.

We start with a well-known fact.

\begin{lemm} \label{lm:p_in_W(Z)}
In the Witt ring $W(\ZZ)$ we have
$$
p=(p,1-p^{p-1},\epsilon_2p^{p-1},\epsilon_3p^{p-1},\epsilon_4p^{p-1},\dots)
$$
where $\epsilon_2,\epsilon_3,\epsilon_4,\dots$ are principal $p$-adic units,
if $p\ge 3$, and
$$
2=(2,-1,\epsilon_2 2^2,\epsilon_32^3,\epsilon_42^5,\dots,\epsilon_n2^{2^{n-2}+1},\dots)
$$
where $\epsilon_2,\epsilon_3,\epsilon_4,\dots$ are $2$-adic units, if $p=2$.
\end{lemm}

\begin{proo}
We start by proving that for $i\ge 1$ we have:
$$
(1-p^{p-1})^{p^i}=
\left\{
\begin{array}{ll}
1-p^{i+p-1}+\frac{p^i-1}{2}p^{i+2(p-1)}+\dots & \mbox{if } p>2, \\
1 & \mbox{if } p=2.
\end{array}
\right.
$$
If $p=2$ this is clear and we assume $p\ge 3$. Now
$(1-p^{p-1})^{p^i}=1-x_1+x_2-\dots+(-1)^{p^i}x_{p^i}$
where $x_j={p^i \choose j}p^{j(p-1)}$ has valuation $v(x_j)=i-v(j)+j(p-1)$
whenever $1\le j\le p^i$.
Let us write $j=up^a$ with $u\ge 1$ prime to $p$ and $a\ge 0$. Then
$v(x_j)=i-a+up^a(p-1)$ which is increasing both as a function of $a$
and as a function of $u$. If $j\ge 2$, then either $u\ge 2$ or
$a\ge 1$. In the first case we have $v(x_j)\ge i+2(p-1)$ and we have
equality for $j=2$. In the second case
we have $v(x_j)\ge i-1+p(p-1)> i+2(p-1)$ since $p\ge 3$. The claim follows.
Now we come to the statement of the lemma itself. The proof for $p=2$ is
similar and we focus on the case $p\ge 3$. The Witt vector
$p=(a_0,a_1,a_2,\dots)$ is determined by the equalities
$$
a_0^{p^n}+pa_1^{p^{n-1}}+\dots+p^{n-1}a_{n-1}^p+p^na_n=p
$$
for all $n\ge 0$. In particular $a_0=p$ and $a_1=1-p^{p-1}$.
By the computation of the $p$-adic first terms of $(1-p^{p-1})^{p^i}$
which we started with, if $n\ge 2$ we have
$$
\frac{p-pa_1^{p^{n-1}}}{p^n}=\frac{p^{n-1+p-1}+\dots}{p^{n-1}}=p^{p-1}+\dots
$$
For $n\ge 2$, by induction using the equality
$$
a_n=\frac{p-pa_1^{p^{n-1}}}{p^n}
-p^{-n}(a_0^{p^n}+p^2a_2^{p^{n-2}}+\dots+p^{n-1}a_{n-1}^p),
$$
we see that the $p$-adic leading term of $a_n$ is $p^{p-1}$.
\end{proo}

\begin{coro} \label{co:ring_O}
Let $\cO=\ZZ[C,\Lambda]/(p-C\Lambda^{p-1})$ and let
$c,\lambda\in\cO$ be the images of $C,\Lambda$. There exists a unique
$d=(d_0,d_1,d_2,\dots)=(c,1-p^{p-1},d_2,\dots)$ in $W^{\lambda^p}(\cO)$ such that
$$
p[\lambda]={\lambda^p}(d)=(\lambda^p d_0,\lambda^pd_1,\lambda^pd_2,\dots).
$$
\end{coro}

\begin{proo}
From the lemma we deduce
$p[\lambda]
=(c\lambda^p,(1-p^{p-1})\lambda^p,\epsilon_2p^{p-1}\lambda^{p^2},
\epsilon_3p^{p-1}\lambda^{p^3},\dots)$.
The coefficients of this vector are divisible by $\lambda^p$, thus
$d_0,d_1,d_2,\dots$ exist. They are unique since $\lambda$ is not a
zero divisor in $\cO$.
\end{proo}

Thus for any $\ZZ$-algebra $R'$ and any elements $c',\lambda'\in R$
satisfying $p=c'\lambda'^{p-1}$ there is a well-determined
$d'\in W^{\lambda'^p}(R')$ such that
$p[\lambda']=(\lambda'^p d'_0,\lambda'^pd'_1,\lambda'^pd'_2,\dots)$.
In particular, our choice of elements $\lambda_i\in R$ satisfying
Assumption~\ref{as:reinf} determines elements
$d_i=(d_{i0},d_{i1},\dots)\in W^{\lambda_i^p}(R)$ such
that
$$
p[\lambda_i]=(\lambda_i^p d_{i0},\lambda_i^pd_{i1},\lambda_i^pd_{i2},\dots).
$$
These are the elements denoted $p\tilde{\lambda}_i/\lambda_i^p$
in \cite{SS2} and $p[\lambda_i]/\lambda_i^p$ in \cite{MRT}.

\begin{noth} {\bf Description of the procedure.}
As in~\ref{ss:fil_gs}, we fix positive integers $L_i$ such that
$(\lambda_i)^{L_i}\in \lambda_{i+1}R$ for all $i\ge 1$ and positive
integers $M,N$. The $n$-th step of the induction produces data:
\begin{itemize}
\item $\bbf{h}^n=(\bbf{a}^n,\bbf{b}^n,\bbf{u}^n,\bbf{v}^n,\bbf{z}^n)$ is
a big frame including two frames of definition
$\bbf{e}^n=(\bbf{a}^n,\bbf{b}^n)$ and $\bbf{f}^n=(\bbf{u}^n,\bbf{v}^n)$
of filtered group schemes and a compatibility between them given by $\bbf{z}^n$,
\item $(\bbf{e}^n,D_{n-1}^{-1},D_{n-1},\cE_n,U^n)$
is a framed group scheme of type $(\lambda_1,\dots,\lambda_{n+1})$,
\item $(\bbf{f}^n,E_{n-1}^{-1},E_{n-1},\cF_n,\bar{U}{}^n)$
is a framed group scheme of type $(\lambda_1^p,\dots,\lambda_{n+1}^p)$,
\item $\Psi^n:\cE_n\to\cF_n$ is an isogeny commuting with the morphism
$\Theta^n$,
\item $\Upsilon^n:(W_{M,N,\lambda_{n+1}})^n\to (W_{M,N,\lambda_{n+1}})^n$
is a matrix of operators (made precise below) describing $(\Psi^n)^*$.
\end{itemize}
The condition that $\Psi^n$ commutes with $\Theta^n$ involves
implicitly the maps $\alpha_{\cE_n}:\cE_n\to (\GG_m)^n$ and
$\beta_{\cF_n}:\cF_n\to (\GG_m)^n$ provided by the construction of framed
group schemes, and may be pictured by the commutative diagramme:
$$
\xymatrix{
\cE_n \ar[r]^{\Psi^n} \ar[d]_{\alpha_{\cE_n}} & \cF_n \ar[d]^{\beta_{\cF_n}} \\
(\GG_m)^n \ar[r]^{\Theta^n} & (\GG_m)^n . \\}
$$
Since $\beta_{\cF_n}$ is an isomorphism on the generic fibre, there is in
any case a rational map $\cE_n\dasharrow \cF_n$. The morphism $\Psi^n$
is determined as the unique morphism extending this rational map.
In fact, the choice of the big frame $\bbf{h}^n$ will guarantee that $\Psi^n$
exists and we may as well remove it from the list above; we included it for
clarity of the picture.
\end{noth}

\begin{noth} {\bf Initialization.}
We set $W^0=(W_{M,N})^0=0$, $\cE_0=\cF_0=0$ and
\begin{itemize}
\item $\bbf{h}^1=(0,0,0,0,0)$,
\item $D_0=D_0^{-1}=1$, $\cE_1=\cG^{\lambda_1}$,
\item $E_0=E_0^{-1}=1$, $\cF_1=\cG^{\lambda_1^p}$,
\item $U^1=F^{\lambda_1}:W_{M,N,\lambda_2}\to W_{M_1,N_1,\lambda_2}$,
\item $\bar U^1=F^{\lambda_1^p}:W_{M,N,\lambda_2^p}
\to W_{M_1,N_1,\lambda_2^p}$,
\item $\Upsilon^1=T_{d_1}:W_{M,N,\lambda_2^p}\to W_{M,N,\lambda_2}$,
\end{itemize}
where $M_1,N_1$ are suitable integers whose existence comes
from Lemmas~\ref{lm:ideal} and \ref{lm:F_stable}.
\end{noth}

\begin{noth} {\bf Induction.}
As usual, we assume that objects in dimension $i$ have been constructed
for $1\le i\le n$ and we explain how to produce
$\bbf{h}^{n+1}=(\bbf{a}^{n+1},\bbf{b}^{n+1},\bbf{u}^{n+1},
\bbf{v}^{n+1},\bbf{z}^{n+1})$
and the related data.

\medskip

\no {\bf A.}
In order to define the big scheme of frames, first we introduce an
$n+1$-dimensional vector $\bbf{c}^{n+1}=(\bbf{a}^n,[\lambda_{n+1}])$.
We recall that Assumption~\ref{as:reinf} is supposed to be satisfied.
The fundamental ingredient of the induction is given by the following result.

\begin{theo}
Let $\cE_{n+1},\cF_{n+1}$ be framed group schemes of types
$(\lambda_1,\dots,\lambda_{n+1})$, $(\lambda_1^p,\dots,\lambda_{n+1}^p)$.
Let $(\bbf{a}^{n+1},\bbf{b}^{n+1})$ and $(\bbf{u}^{n+1},\bbf{v}^{n+1})$
be the defining frames. Assume that the Kummer subscheme $G_n\subset \cE_n$
is finite locally free and that the rational map $\cE_n\dasharrow\cF_n$
extends to an isogeny with kernel $G_n$. Then, the
following conditions are equivalent:
\begin{trivlist}
\itemn{1} the Kummer subscheme $G_{n+1}\subset \cE_{n+1}$ is finite
locally free and the rational map $\cE_{n+1}\dasharrow\cF_{n+1}$ extends
to an isogeny with kernel $G_{n+1}$,
\itemn{2} there exists $\bbf{z}^{n+1}\in (W^{\lambda_{n+1}^p})^n(R)$
such that \
$p\bbf{a}^{n+1}-\bbf{c}^{n+1}-\Upsilon^{n}\bbf{u}^{n+1}
={\lambda_{n+1}^p}(\bbf{z}^{n+1})$.
\end{trivlist}
\end{theo}

\begin{proo}
This is proven in \cite{MRT}, Theorem 7.1.1,
in the case where the ring $R$ is a discrete
valuation ring, with uniformizer $\pi$. The proof uses general power
series computations and it is clear while reading it that it works for
an arbitrary $\ZZ_{(p)}$-algebra $R$ satisfying our assumptions. We indicate
the necessary changes of notation: $\lambda_i$ has to be replaced by
$\pi^{l_i}$, and ${\lambda_{n+1}^p}(\bbf{z}^{n+1}_i)$ has to be
replaced by $T_{\bbf{z}_i^{n+1}}([\pi^{pl_{n+1}}])$, where
$\bbf{z}^{n+1}=(\bbf{z}^{n+1}_1,\dots,\bbf{z}^{n+1}_n)$.
\end{proo}

Given this theorem, we can choose a big frame
$$
\bbf{h}^{n+1}=(\bbf{a}^{n+1},\bbf{b}^{n+1},\bbf{u}^{n+1},
\bbf{v}^{n+1},\bbf{z}^{n+1})
$$
living in a big scheme of frames whose heavy but obvious definition we omit.

\medskip

\no {\bf B.}
Using the components $\bbf{a}^{n+1}$ and $\bbf{u}^{n+1}$ of the frame,
we define:
$$
\begin{array}{lcl}
D_n & =
& \prod_{i=1}^n\,E^{L_i,M_i,N_i}_p(\bbf{a}^{n+1}_i,\lambda_i,D_{i-1}^{-1}X_i)
\bigskip \\
D_n^{-1}& =
& \prod_{i=1}^n\,E^{L_i,M_i,N_i}_p(-\bbf{a}^{n+1}_i,\lambda_i,D_{i-1}^{-1}X_i)
\bigskip \\
E_n & =
& \prod_{i=1}^n\,E^{L_i,M_i,N_i}_p(\bbf{u}^{n+1}_i,\lambda_i^p,E_{i-1}^{-1}Y_i)
\bigskip \\
E_n^{-1}& =
& \prod_{i=1}^n\,E^{L_i,M_i,N_i}_p(-\bbf{u}^{n+1}_i,\lambda_i^p,E_{i-1}^{-1}Y_i). \\
\end{array}
$$

\medskip

\no {\bf C.} At this step, we define $\cE_{n+1}$ and $\cF_{n+1}$
in the same way as in~\ref{induct}, Step C.

\medskip

\no {\bf D.} At this step, we define morphisms
$$
U^{n+1},\bar{U}^{n+1}:(W_{M,N,\lambda_{n+2}})^{n+1}
\to (W_{M_{n+1},N_{n+1},\lambda_{n+2}})^{n+1}
$$
like in~\ref{induct}, Step D, and the operator
$\Upsilon^{n+1}:(W_{M,N,\lambda_{n+2}})^{n+1}\to (W_{M,N,\lambda_{n+2}})^{n+1}$
by the matrix
$$
\Upsilon^{n+1}=\left(
\begin{array}{cccc}
& & & -T_{\bbf{z}^{n+1}_1} \\
& \Upsilon^n & & \vdots \\
& & & -T_{\bbf{z}^{n+1}_n} \\
0 & \dots & 0 & T_{d_{n+1}} \\
\end{array}
\right).
$$
This concludes the inductive construction.
\end{noth}

\begin{theo} \label{th:universal_Kummer}
Let $\BB^n_\nu$ be the finite flat covers of affine space $\AA^n$
defined in~\ref{th:universal_alg}.
There exists a sequence indexed by $\nu\ge 1$ of affine
$\BB^n_\nu$-schemes $\sK_n^\nu=\Spec(\sM_n^\nu)$ of finite type,
without $\Lambda$-torsion, framed $\sM_n^\nu$-group schemes
$\sE_n^\nu$ of type $(\Lambda_1,\dots,\Lambda_n)$ and
$\sF_n^\nu$ of type $(\Lambda_1^p,\dots,\Lambda_n^p)$, and an
isogeny $\sE_n^\nu\to\sF_n^\nu$ with finite locally free kernel $\sG_n^\nu$
compatible with the maps to $(\GG_m)^n$. This isogeny is universal
in the same sense as in~\ref{th:universal_alg}.
\end{theo}

\begin{proo}
Omitted.
\end{proo}

The family $\sG_n=(\sG_n^\nu)_{\nu\ge 1}$ is a finite flat group scheme
over the ind-scheme $(\sK_n^\nu)_{\nu\ge 1}$. We call it the
{\em universal Kummer group scheme}.

\medskip

We conclude with a remark on the operator $\Upsilon^n$. By construction,
it represents the pullback $(\Psi^n)^*$, which implies that modulo
$\lambda_{n+1}$ it maps the
subspace $\ker(\bar{U}{}^n)$ into the subspace $\ker(U^n)$. In fact, we can
do better: it is possible to include in the induction the construction of
a matrix $\Omega^n$ such that $U^n\,\Upsilon^n=\Omega^n\,\bar{U}{}^n$. This
is a reflection of the fact that among the morphisms from a filtered group
to $\GG_m$, not only the {\em group morphisms} (represented by $\ker(U^n)$)
but also the {\em fundamental morphisms} (represented by the ambient $W^n$)
are meaningful. On the diagonal, the entries of the matrix $\Omega^n$ should
be operators $T'_{d_i}$ (see below) satisfying
$F^{\lambda_i}\circ T_{d_i}=T'_{d_i}\circ F^{\lambda_i^p}$.
In fact, these matrices are defined by $\Omega^1=T'_{d_1}$ and
$$
\Omega^{n+1}=\left(
\begin{array}{cccc}
& & & * \\
& \Omega^n & & \vdots \\
& & & * \\
0 & \dots & 0 & T'_{d_{n+1}} \\
\end{array}
\right).
$$
We do not want to go into the full details of the construction of $\Omega^n$.
We simply note that the essential task is to define the diagonal entries
$T'_{d_i}$. We end the paper with the proof of existence and unicity of
these endomorphisms.

\begin{lemm}
Let $\cO=\ZZ[C,\Lambda]/(p-C\Lambda^{p-1})$ and let $c,\lambda\in\cO$ be the
images of $C,\Lambda$. Let $d=(c,1-p^{p-1},\dots)$ be the unique vector such
that $p[\lambda]=(\lambda^p d_0,\lambda^pd_1,\dots)$, as in
Corollary~\ref{co:ring_O}. Then there exists a unique endomorphism $T'_d:W\to W$
such that $F^\lambda\circ T_d=T'_d\circ F^{\lambda^p}$ as endomorphisms of
the $\cO$-group scheme $W$.
\end{lemm}


\begin{proo}
Since $F^{\lambda^p}$ is an epimorphism, then $T'_d$ is unique and
we only have to prove that it exists. Let $\Phi:W\to (\GG_a)^{\NN}$ be the
Witt morphism of $\cO$-ring schemes. Given that the schemes $\Spec(\cO)$ and $W$
have no $p$-torsion, the morphism $\Phi$ is a monomorphism and it is enough to look
for $T'_d:W\to W$ such that
$\Phi\circ F^\lambda\circ T_d=\Phi\circ T'_d\circ F^{\lambda^p}$.
Let $f$ and $t_d$ be the endomorphisms of $(\GG_a)^{\NN}$ such that
$\Phi\circ F=f\circ\Phi$ and $\Phi\circ T_d=t_d\circ\Phi$. They are defined by:
\begin{itemize}
\item $f(x_0,x_1,x_2,\dots)=(x_1,x_2,x_3,\dots)$,
\item $t_d(x_0,x_1,x_2,\dots)=(y_0,y_1,y_2,\dots)$ with
$y_n=d_0^{p^n}x_n+pd_1^{p^{n-1}}x_{n-1}+\dots+p^nd_nx_0$.
\end{itemize}
We first construct $t'_d:(\GG_a)^{\NN}\to (\GG_a)^{\NN}$ such that
$$
\big(f-\Phi([\lambda^{p-1}])\Id\big)\circ t_d
=t'_d\circ \big(f-\Phi([\lambda^{(p-1)p}])\Id\big).
$$
Let $y=(y_0,y_1,y_2,\dots)$ be a Witt vector of indeterminates and
let us write
$$
\begin{array}{r}
\left(\left[f-\Phi([\lambda^{p-1}])\right]\circ t_d\right)(y)
= (\alpha_0,\alpha_1,\alpha_2,\dots), \bigskip \\
\left[f-\Phi([\lambda^{(p-1)p}])\right](y)
= (\beta_0,\beta_1,\beta_2,\dots).
\end{array}
$$
Given that $\Phi([a])=(a,a^p,a^{p^2},\dots)$ and that addition and
multiplication in $(\GG_a)^{\NN}$ are componentwise, we compute:
\begin{align*}
\alpha_n = & \big(d_0^{p^{n+1}}y_{n+1}+pd_1^{p^n}y_n+\dots
+p^nd_n^py_1+p^{n+1}d_{n+1}y_0\big) \\
& -\lambda^{p^n(p-1)}
\big(d_0^{p^n}y_n+pd_1^{p^{n-1}}y_{n-1}+\dots+p^{n-1}d_{n-1}^py_1+p^nd_ny_0\big)
\end{align*}
and $\beta_n=y_{n+1}-\lambda^{p^{n+1}(p-1)}y_n$ for all $n\ge 0$. The existence
of $t'_d$ means that $\alpha_n$ is a polynomial with coefficients in $\cO$
in the variables $\beta_0,\beta_1,\beta_2,\dots$ for each $n$. Since the
$\alpha_n$ and $\beta_n$ are linear in $y$, this in turn means that
we get $\alpha_n=0$ under the specializations
$$
y_1=\lambda^{p(p-1)}y_0 \ , \
y_2=\lambda^{p^2(p-1)}y_1 \ ,\dots, \
y_{i+1}=\lambda^{p^{i+1}(p-1)}y_i \ ,\dots
$$
This amounts to $y_i=\lambda^{p(p^i-1)}y_0$ for each $i$. Now
$$
\alpha_n\left(y_0,\lambda^{p(p-1)}y_0,\lambda^{p(p^2-1)}y_0,
\lambda^{p(p^3-1)}y_0,\dots\right)
$$
is equal to $y_0$ times
\begin{align*}
& \big(d_0^{p^{n+1}}\lambda^{p(p^{n+1}-1)}+pd_1^{p^n}\lambda^{p(p^n-1)}
+\dots+p^nd_n^p\lambda^{p(p-1)}+p^{n+1}d_{n+1}\big) \\
& -\lambda^{p^n(p-1)}
\big(d_0^{p^n}\lambda^{p(p^n-1)}+pd_1^{p^{n-1}}\lambda^{p(p^{n-1}-1)}
+\dots+p^{n-1}d_{n-1}^p\lambda^{p(p-1)}+p^nd_n\big).
\end{align*}
If we recall that $p\lambda^{p^i}=
\lambda^{p^{i+1}}d_0^{p^i}+p\lambda^{p^i}d_1^{p^{i-1}}+\dots+p^i\lambda^pd_i$
for all $i$ by definition of $d$, then we indeed find that this quantity
vanishes. This proves the existence of $t'_d$ as required.
In order to find a morphism $T'_d$ such that $\Phi\circ T'_d=t'_d\circ\Phi$,
we use Bourbaki \cite{B}, \S~1, no.~2, Prop.~2, applied to $t'_d\circ\Phi$,
viewed as a sequence of elements in the ring $H^0(W,\cO_W)=\cO[Z_0,Z_1,\dots]$
endowed with the endomorphism raising each variable to the $p$-th power.
\end{proo}

\appendix

\section{Appendix: errata for the paper~\cite{SS2}}
\label{app:errata}

This is a list of typographical slips that we are aware of in the preprint
{\em On the unified Kummer-Artin-Schreier-Witt theory} \cite{SS2} by
T. Sekiguchi and N. Suwa.
We thank heartily Guillaume
Pagot and Michel Matignon who communicated to us their personal list of errata,
which is included here. Notation "p.~$x$ l.~$y$" means page $x$, line $y$
from the top (not counting running headers) and "p.~$x$ l.~$-y$" means
page $x$, line $y$ from the bottom (counting equations as one line).

\bigskip

\tiret p.~4, l.~-7: replace $\lambda_n$ by $\lambda$.

\tiret p.~4, l.~-1: replace $\lambda_n$ by $\lambda_{n+1}$, two times.

\tiret p.~5, l.~3: replace $\lambda_n$ by $\lambda_{n+1}$, two times.

\tiret p.~13, l.~8: slight conflict of notation between
$\alpha^{(\ell)}:\cE_\ell\to (\GG_{m,A})^\ell$
and $\alpha^{(\lambda)}:\cG^{(\lambda)}\to\GG_m$.

\tiret p.~14, Lemma 3.4: replace $\frac{X_i}{D_i(\YY)}$ by
$\frac{Y_i}{D_i(\XX)}$
two times in the statement, three times in the proof.

\tiret p.~18, l.~-7: replace $W_n(A)$ by $W(A)$.

\tiret in various places, functors (or morphisms between functors) are defined
by one of the following two procedures: define the functor on elements
$a\in A$ for varying coefficient rings $A$, or define the functor in the
universal case. These procedures are of course equivalent, but the reader
should be a little careful because several times the definition is given in
terms of one procedure while the actual occurrence of the functor in the text
is in terms of the other procedure. Here are some examples:
\begin{trivlist}
\itemn{i}
p.~19 the morphism $[p]:W\to W$ is defined for a vector $\bbf{b}\in W(A)$
but its first apparition on p.~25 is $[p]\UU$ for some indeterminate vector $\UU$.
\itemn{ii}
p.~20 the map $\bbf{x}\mapsto T_{\bbf{a}}\bbf{x}$ is defined on p.~20 for
$\bbf{a},\bbf{x}\in W(A)$ but later on p.~26 there is $T_{\VV}\WW$ for
indeterminate vectors $\VV,\WW$.
\itemn{iii}
p.~20 the notation $\bbf{a}/\lambda$ is introduced for $\bbf{a}\in W(A)$
but on p.~25 there is $\frac{1}{\Lambda_2}\UU$.
\itemn{iv}
p.~23 the deformed exponential $E_p(\bbf{a},\lambda,X)$ is defined for
$\bbf{a}\in W(A)$ and $\lambda\in A$, but on p.~25 there is
$E_p(\UU,\Lambda,X)$.
\end{trivlist}

 \tiret p.~22, Lemma 4.7: read $\sigma(f)\equiv f^p
\mod pA$.

\tiret p.~23 the deformed exponential $E_p(\bbf{a},\lambda,X)$ is
defined for $\bbf{a}\in W(A)$ and $\lambda\in A$, but on p.~25
there is $E_p(\UU,\Lambda,X)$.

\tiret p.~25, Lemma 4.10: the symbol $\UU^{(p^k)}$ is defined on page 8.

\tiret p.~26, l.~11: the second factor should be
$
\left(\prod_{r\ge 1}\left(1+\Lambda_1^{p^r}X^{p^r}\right)
^{\frac{1}{p^r\Lambda_1^{p^r}}
\Phi_{r-1}(F^{(\Lambda_1)}\UU)}\right)^{\frac{W_0}{\Lambda_2}}$.

 \tiret p.~26 l.~-7: replace
$\Lambda_1^{p^{r+k}}X^{p^{r+k}}$ by $\Lambda_1^pX^p$.


\tiret p.~27 l.~6: replace $\frac{1}{\Lambda_2^{\ell-1}}$ by
$\frac{1}{\Lambda_2^{p^{\ell-1}}}$.

\tiret p.~30 l.~-7: replace $d_1$ by $d_0$.



\tiret p.~37 l.~3: replace "definition (51)" by "definition (45)".

\tiret p.~40 l.~-7: replace $\YY=\YY$ by $\YY^1=\YY$.

\tiret p.~42 l.~11: read
$\bbf{b}_1^2=\frac{1}{\lambda_2}F^{(\lambda_1)}\bbf{a}^1_1$.

\tiret p.~43, definition of $D_{k+1}$: in order to define the group scheme
$\cE_n$ one needs polynomials, but the $D_k$ as defined are power series
since the Witt vectors $\bbf{a}_i^{k+1}$ have coefficients in $A$. In fact,
the good way to proceed is to first consider
$$
\prod_{i=1}^{k+1}\,
E_p\left(\bbf{a}_i^{k+1},\lambda_i,\frac{X_{i-1}}{D_{i-1}(X_0,\dots,X_{i-1})}\right)
$$
as a polynomial with coefficients in $A/\lambda_{k+2}$ and then lift it to a
polynomial $D_{k+1}$ of the same degree with coefficients in $A$.


\tiret p.~46 l.~1: replace $\XX$ by $X$. Note that of course $U^1$ is relative
to $\lambda_n$, not $\lambda_1$.


\tiret p.~47, lines 3 to 6: Proposition~3.5 requires $B$ to be a discrete
valuation ring.

\tiret p.~50, in the statement and proof of Lemma~6.2: replace
$M^{p^r}$ by $M^{p^{r+\ell}}$.

\tiret p.~52 l.~5: replace $X_{i-1}+Y_{i-1}$ by $X_{k-1}+Y_{k-1}$.

\tiret p.~59 l.~-4: replace
$D_1(c_0)+\lambda_{(1)}=\zeta_2$ by $D_1(c_0)+\lambda_{(1)}c_1=\zeta_2$.

\tiret p.~62 l.~5: replace $F^n\bar{\bbf{a}}^n_1$ by
$F^{(\lambda)}\bar{\bbf{a}}^n_n$.

\tiret p.~64 l.~3: replace $L_p(U)$ by $L_{p,N}(U)$.

\tiret p.~64 l.~-4: replace
$F^{(\lambda)}\bar{\bbf{a}}^n_{n-1}$ by $F^{(\lambda)}\bar{\bbf{a}}^n_n$.

\tiret p.~67: there are some incompatibilities between $\alpha_k$ in
Definition 8.1 and $\alpha_1,\alpha_2$ at the beginning of p.~67.
Same thing on p.~80.

\tiret p.~70 l.~4: In the formula (89) replace
$\frac{1}{\bar{D}_1(Y_0)+\lb^p Y_{1}}$ by $\frac{1}{1+\lb^p Y_0}$.
And moreover the $\bar{D}_n$ should be a polynomial and not a
series: see comment (17).

 \tiret p.~71 l.~7: replace $t_2^pt_2^{-1}$ by
$t_2^pt_1^{-1}$.

\tiret p.~72, l.~3: replace $A_\lambda^p$ by $A_{\lambda^p}$, two times.

\tiret p.~74 l.~1: replace $\End(W(A))$ by $\End(W(A)^n)$.

\tiret p.~74 l.~8: replace $W(A)^{n+1}$ by $W(A/\lambda^p)^{n+1}$.

\tiret p.~74 l.~11: in the second component of the vector
$T_{(p\bbf{a}^n-\bbf{c}^{n-1}-\Upsilon^n\bbf{u}^n)/\lambda^p}$, replace
$\bbf{a}_1^n$ by $\bbf{a}_2^n$.

\tiret p.~74 l.~-1: replace $K_n^{W_0/N}$ by $\widetilde{K_n}^{W_0/N}$.

\tiret p.~74 l.~-1 and p.~75 l.~2: $\tilde E_p$ was defined
in (32), p.~26 and $G_p$ was defined in (34), p.~29.

\tiret p.~75: in the statement of Theorem~9.4, the target of the map
$\Upsilon^n$ is $\ker(U^n)$.

\tiret p.~83 l.~6: read
$\equiv (p\alpha_2-\alpha_1,(1-p^{p-1}\alpha_2^p)-\gamma+C_1(p\alpha_2,-\alpha_1))$.

\tiret p.~83 l.~9:
read $\delta_1:=(1-p^{p-1})\alpha_2^p-\gamma+C_1(p\alpha_2,-\alpha_1)$.

\tiret p.~83 l.~10: replace $v(\alpha^2)$ by $=v(\alpha_2^p)$.







\end{document}